\documentclass[11pt,reqno,twoside]{article}
\def\COMPLETE{}
\input{package_header} 
\begin{document}
\title{Inexact and Stochastic Generalized Conditional Gradient with Augmented Lagrangian and Proximal Step}
\author{Antonio Silveti-Falls\thanks{Normandie Universit\'e, ENSICAEN, UNICAEN, CNRS, GREYC, France. E-mail: tonys.falls@gmail.com, cecio.molinari@gmail.com, Jalal.Fadili@ensicaen.fr.} \and
Cesare Molinari\samethanks \and 
Jalal Fadili\samethanks
}
\date{}
\maketitle
\begin{flushleft}\end{flushleft}
\begin{abstract}
In this paper we propose and analyze inexact and stochastic versions of the CGALP algorithm developed in \cite{silveti}, which we denote \icgalp, that allows for errors in the computation of several important quantities. In particular this allows one to compute some gradients, proximal terms, and/or linear minimization oracles in an inexact fashion that facilitates the practical application of the algorithm to computationally intensive settings, e.g. in high (or possibly infinite) dimensional Hilbert spaces commonly found in machine learning problems. The algorithm is able to solve composite minimization problems involving the sum of three convex proper lower-semicontinuous functions subject to an affine constraint of the form $Ax=b$ for some bounded linear operator $A$. Only one of the functions in the objective is assumed to be differentiable, the other two are assumed to have an accessible prox operator and a linear minimization oracle. As main results, we show convergence of the Lagrangian to an optimum and asymptotic feasibility of the affine constraint as well as weak convergence of the dual variable to a solution of the dual problem, all in an almost sure sense. Almost sure convergence rates, both pointwise and ergodic, are given for the Lagrangian values and the feasibility gap. Numerical experiments verifying the predicted rates of convergence are shown as well.
\end{abstract}




\ifdefined\COMPLETE
\else
\documentclass[12pt]{article}
\input{package_header} 
\begin{document}
\fi
\section{Introduction}
\begin{subsection}{Problem Statement}
We consider the following composite minimization problem,
\begin{equation}\label{PProb}\tag{$\mathrsfs{P}$}
\min\limits_{x\in\HH_p}\brac{f\para{x} + g\para{Tx} + h\para{x}: Ax=b},
\end{equation}
and its associated \emph{dual problem},
\begin{equation}\label{DProb}\tag{$\mathrsfs{D}$}
\min\limits_{\mu\in\HH_d}\para{f+g\circ T+h}^*\para{-A^*\mu}+\ip{\mu,b}{},
\end{equation}
where we have denoted by $*$ both the \emph{Legendre-Fenchel conjugate} and the \emph{adjoint operator}, to be understood from context.
We consider $\HH_p$, $\HH_d$, and $\HH_v$ to be arbitrary real Hilbert spaces, possibly infinite-dimensional, whose indices correspond to a primal, dual, and auxilliary space, respectively; $A:\HH_p\to\HH_d$ and $T:\HH_p\to\HH_v$ to be bounded linear operators with $b\in\ran(A)$; functions $f$, $g$, and $h$ to all be convex, closed, and proper real-valued functions. Additionally, we will assume that the function $f$ satisfies a certain differentiability condition generalizing Lipschitz-smoothness, H\"{o}lder-smoothness, etc (see \defref{def_smooth}), that the function $g$ has a proximal mapping which is accessible, and that the function $h$ admits an accessible linearly-perturbed minimization oracle with $C\eqdef \dom\para{h}$ a weakly compact subset of $\HH_p$. 

In fact, the problem under consideration here is exactly the same as that of \cite{silveti}, however, in this work, we consider an inexact extension of the algorithm presented and analyzed in \cite{silveti} to solve \eqref{PProb}. The extension amounts to allowing either deterministic or stochastic errors in the computation of several quantities, including the gradient or prox terms, e.g. $\nabla f$, $\prox_{\beta g}$, and the linear minimization oracle itself.
\end{subsection}

\begin{subsection}{Contribution and prior work}
The primary contribution of this work is to analyze inexact and stochastic variants of the \cgal algorithm presented in \cite{silveti} to address \eqref{PProb}. We coin this algorithm \textbf{I}nexact \textbf{C}onditional \textbf{G}radient with \textbf{A}ugemented \textbf{L}agrangian and \textbf{P}roximal-step (\icgalp). Although there has been a great deal of work on developing and analyzing Frank-Wolfe or conditional gradient style algorithms in both the stochastic and deterministic case, e.g. \cite{Hazan12,Hazan16,sfw3,sfw4,sfw5,sfw6,sfw1,sfw2}, or \cite{sfw7}, little to no work has been done to analyze the generalized version of these algorithms for nonsmooth problems or problems involving an affine constraint, as we consider here. To the best of our knowledge, the only such work is \cite{loca}, where the authors consider a stochastic conditional gradient algorithm applied to a composite problem allowing nonsmooth terms. The nonsmooth term is possibly an affine constraint but it is addressed through smoothing rather than through an augmented Lagrangian with a dual variable, in contrast to our work.

We show asymptotic feasibility of the primal iterates for the affine constraint, convergence of the Lagrangian values at each iteration to an optimum value, weak convergence of the sequence of dual iterates to a solution of the dual problem, and provide worst-case rates of convergence for the feasibility gap and the Lagranian values. The rates of convergence are given both subsequentially in the pointwise sense and globally, i.e. for the entire sequence of iterates, in the ergodic sense where the Ces\'aro means are taken with respect to the primal step size. In the case where \eqref{PProb} admits a unique solution, we furthermore have that the sequence of primal iterates converges weakly to the solution. These results are shown to hold almost surely and are established for a family of parameters satisfying abstract open loop conditions, i.e. sequences of parameters which do not depend on the iterates themselves. We exemplify the framework on problem instances involving a smooth risk minimization where the gradient is computed inexactly either with stochastic noise or a deterministic error. In the stochastic case, we show that our conditions outlined in \secref{sec:alg} for convergence are satisfied via increasing batch size or variance reduction. In the deterministc setting for minimizing an empirical risk, a sweeping approach is described.
\end{subsection}

\begin{subsection}{Organization}
The remainder of the paper is divided into four sections. In \secref{sec:not} the necessary notation and prior results are recalled, consisting primarily of convex analysis, real analysis, and elementary probability. In \secref{sec:alg} the assumptions on the problem structure and the parameters are noted, the \icgalp algorithm itself is presented. In \secref{sec:mainres}, the main results, e.g. feasibility, Lagrangian convergence, and rates, are established. The analysis and results extend those of \cite{silveti} to the inexact and stochastic setting. In \secref{sec:exam} and \secref{sec:ex2}, we consider different problem instances where inexact deterministic or stochastic computations are involved. Numerical results are reported in \secref{sec:numexp} to support our theoretical findings. Finally, in \secref{sec:conc}, we summarize the work and provide some closing remarks.
\end{subsection}

\ifdefined\COMPLETE
\else
\end{document}
\fi

\ifdefined\COMPLETE
\else
\documentclass[12pt]{article}
\input{package_header} 
\begin{document}
\fi
\section{Notation and Preliminaries}\label{sec:not}
Many of the following notations for probabilistic concepts are adopted from \cite{combstoch}. We denote by $\prspace$ a probability space with set of events $\events$, $\sigma$-algebra $\sigalg$, and probability measure $\prob$. When discussing random variables we will assume that any Hilbert space $\mathcal{H}$ is endowed with the Borel $\sigma$-algebra, $\borel{\mathcal{H}}$. We denote a \emph{filtration} by $\Filt = \seq{\filt_k}$, i.e. a sequence of sub-$\sigma$-algebras which satisfies $\filt_k\subset \filt_{k+1}$ for all $k\in\mathbb{N}$. Given a set of random variables $\brac{a_0,\ldots,a_n}$, we denote by $\sigma\para{a_0,\ldots,a_n}$ the $\sigma$-algebra generated by $a_0,\ldots,a_n$. An expression $\para{P}$ is said to hold $\Pas$ if $\prob\para{\{\omega\in\events: \para{P}\mbox{ holds}\}}=1$. Throughout the paper, both equalities and inequalities involving random quantities should be understood as holding $\mathbb{P}$-almost surely, whether or not it is explicitly written.
\begin{definition}
Given a filtration $\Filt$, we denote by $\ell_+\para{\Filt}$ the set of sequences of $[0,+\infty[$-valued random variables $\seq{a_k}$ such that, for each $k\in\N$, $a_k$ is $\filt_k$ measurable. Then, we also define the following set,
\newq{
\ell^1_+\para{\Filt} \eqdef \brac{\seq{a_k}\in\ell_+\para{\Filt}:\sum\limits_{k\in\N} a_k<+\infty \Pas} 
}
\end{definition}

\begin{lemma}\label{combstoch}
Given a filtration $\Filt$ and the sequences of random variables $\seq{r_k}\in\ell_+\para{\Filt}$, $\seq{a_k}\in\ell_+\para{\Filt}$, and $\seq{z_k}\in\ell_+^1\para{\Filt}$ satisfying,
\newq{
\EX{r_{k+1}}{\filt_k} - r_k \leq -a_k + z_k \Pas
}
then $\seq{a_k}\in\ell^1_+\para{\Filt}$ and $\seq{r_k}$ converges $\Pas$ to a random variable with value in $[0,+\infty[$.
\end{lemma}
\begin{proof}
See \cite[Theorem 1]{robsieg}.
\end{proof}

\begin{lemma}\label{barty}
Given a filtration $\Filt$ and a sequence of random variables $\seq{w_k}\in\ell_+\para{\Filt}$ and a sequence of real numbers $\seq{\gamma_k}\in\ell_+$ such that $\seq{\gamma_kw_k}\in\ell_+^1\para{\Filt}$ and $\seq{\gamma_k}\not\in\ell^1$, then,
\be[label=(\roman*)]
	\item there exists a subsequence $\subseq{w_{k_j}}$ such that
	\newq{
	w_{k_j}\leq \Gamma_{k_j}^{-1} \Pas
	}
	 where $\Gamma_{k_j} = \sum\limits_{n=1}^{k_j} \gamma_n$. In particular, $\liminf\limits_k w_k =0$ $\Pas$.
	\item Furthermore, if there exists a constant $\alpha>0$ such that $w_k-\EX{w_{k+1}}{\filt_{k}}\leq \alpha \gamma_k$ $\Pas$ for every $k\in\mathbb{N}$, then
	\newq{
	\lim\limits_{k}w_k =0 \Pas.
	}
\ee
\end{lemma}
\begin{proof}
The main result is directly from \cite[Lemma 2.2]{barty} and the rates follow from \cite{Alber} trivially extended to the stochastic setting.
\end{proof}

We denote by $\Gamma_0\para{\HH}$ the set of proper, convex, and lower semi-continuous functions $f:\HH\to\R\cup\brac{+\infty}$. We also consider the \emph{domain} of a function $f$ to be $\dom\para{f} \eqdef \brac{x\in\HH: f\para{x}<+\infty}$ and the \emph{Legendre-Fenchel conjugate} of $f$ to be the function $f^*:\HH\to\R\cup\brac{+\infty}$ such that, $\forall y\in\HH$,
\newq{
f^*\para{y} \eqdef \sup\limits_{x\in\HH}\brac{\ip{y,x}{} -f\para{x}}.
}
The \emph{proximal mapping} (or \emph{proximal operator}) associated to the function $f$ with parameter $\beta$ is given by,
\newq{
\prox_{\beta f} \para{x} \eqdef \argmin\limits_{y\in\HH} \brac{f\para{y}  +\frac{1}{2\beta}\norm{x-y}{}^2}.
}
The following elementary result from convex analysis regarding proximal mappings will be used in the proof of optimality.
\begin{proposition}\label{goodineq}
Let $f\in\Gamma_0\para{\HH}$ and denote $x^+ = \prox_f\para{x}$. Then, for all $y\in\HH$,
\newq{
2\para{f\para{x^+} - f\para{x}} +\norm{x^+-y}{}^2 - \norm{x-y}{}^2 + \norm{x^+-x}{}^2\leq 0.
}
\end{proposition}
\begin{proof}
The result is classical and the proof is readily available, e.g. in \cite[Chapter 6.2.1]{PEYPOU}.
\end{proof}
The \emph{subdifferential} of a function $f$ is the set-valued operator $\partial f: \HH\to 2^\HH$ such that, for every $x\in\HH$,
\nnewq{\label{subdiff}
\p f \para{x} \eqdef \brac{u\in\HH: f\para{y}\geq f\para{x} + \ip{u,y-x}{}\quad \forall y\in\HH}
}
We denote $\dom\para{\p f}\eqdef \brac{x\in\HH: \p f\para{x}\neq\emptyset}$ as the \emph{domain of the subdifferential}. For $x\in\dom\para{\p f}$, the \emph{minimal norm selection} of $\p f\para{x}$ is denoted by $\sbrac{\p f\para{x}}^0\eqdef \argmin\limits_{y\in\p f\para{x}}\norm{y}{}$.
The \emph{Moreau envelope} of the function $f$ with parameter $\beta$ is given by,
\newq{
f^{\beta}\para{x} \eqdef \inf\limits_{y\in\HH} \brac{f\para{y} + \frac{1}{2\beta}\norm{x-y}{}^2}.
}
The following proposition recalls some key properties of the Moreau envelope which we will utilize in the analysis of the algorithm.
\begin{proposition}[Moreau envelope properties]\label{moreau}
Given a function $f\in\Gamma_0\para{\HH}$, the following holds:
\begin{enumerate}[label=(\roman*)]
\item\label{moreauclimconv} The Moreau envelope, $f^\beta$, is convex, real-valued, and continuous.
\item\label{moreauclaim1} Lax-Hopf formula: the Moreau envelope is the viscosity solution to the following Hamilton Jacobi equation:
\begin{equation}
\label{LH}
\begin{cases}\frac{\p}{\p \beta}f^{\beta}\para{x} = -\frac{1}{2}\norm{\nabla_x f^{\beta}\para{x}}{}^2 \quad & \para{x,\beta}\in \HH\times (0,\pinfty)\\ f^0\para{x} = f\para{x}\quad & x\in\HH. \end{cases}
\end{equation}
\item\label{moreauclaim2} The gradient of the Moreau envelope, $\nabla f^\beta$, is $\frac{1}{\beta}$-Lipschitz continuous and is given by the expression
\newq{
	\nabla_x f^{\beta}\para{x} = \frac{x-\prox_{\beta f}\para{x}}{\beta}.
}
\item\label{moreauclaim3} $\forall x \in \dom(\p f)$, $\norm{\nabla f^{\beta}\para{x}}{} \nearrow \norm{\sbrac{\p f\para{x}}^0}{}$ as $\beta\searrow 0$.
\item\label{moreauclaim4} $\forall x \in \HH$, $f^{\beta}(x) \nearrow f(x)$ as $\beta\searrow 0$. In addition, given two positive real numbers $\beta'<\beta$, for all $x\in \HH$ we have
\begin{equation*}
\begin{split}
0 \leq f^{\beta'}\para{x}-f^{\beta}\para{x} & \leq \frac{\beta-\beta'}{2}\norm{\nabla_x f^{\beta'}\para{x}}{}^2;\\
0 \leq f\para{x}-f^{\beta}\para{x} & \leq \frac{\beta}{2}\norm{\sbrac{\p f\para{x}}^0}{}^2.
\end{split}
\end{equation*}
\end{enumerate}
\end{proposition}
Given a closed, convex set $\C$, we write $d_\C \eqdef \sup\limits_{x,y\in\C}\norm{x-y}{}$ to denote the \emph{diameter} of $\C$. We denote the \emph{Bregman divergence} of a differentiable, function $F$ by,
\newq{
D_F\para{x,y} \eqdef F\para{x} - F\para{y} - \ip{\nabla F\para{y}, x-y}{}.
}
\begin{definition}[$\para{F,\zeta}$-smoothness]\label{def_smooth}
Let $F:\HH\to\R\cup\brac{+\infty}$ and $\zeta:]0,1]\to\R_+$. The pair $\para{f,\C}$, where $f:\HH\to\R\cup\brac{+\infty}$ and $\C\subset \dom\para{f}$, is said to be $\para{F,\zeta}$-smooth if there exists an open set $\C_0$ such that $\C\subset\C_0\subset\inte\para{\dom\para{F}}$ and,
\be[label=(\roman*)]
\item $F$ and $f$ are differentiable on $\C_0$;\label{Ffdiff}
\item $F-f$ is convex on $\C_0$;\label{F-fconv}
\item it holds \label{FK}
\nnewq{
\KFC \eqdef \sup\limits_{\sst{x,s\in\C; \gamma\in]0,1]\\ z=x+\gamma\para{s-x}}} \frac{D_F\para{z,x}}{\zeta\para{\gamma}}<+\infty.
}
\ee
\end{definition}
\begin{remark}\label{fremark}
An important consequence of \defref{def_smooth}\ref{Ffdiff} and \defref{def_smooth}\ref{F-fconv} in $\para{F,\zeta}$-smoothness is the following. Let $\para{f,\C}$ be $\para{F,\zeta}$ smooth. Then, for any $x,y\in\C$, we have,
\newq{
f\para{y}\leq f\para{x} + \ip{\nabla f\para{x}, y-x}{} + D_F\para{y,x}.
}
Moreover, by \defref{def_smooth}\ref{FK}, if $y = x +\gamma\para{s-x}$ for some $s\in\C$ and $\gamma\in]0,1]$, we have,
\nnewq{
D_F\para{y,x} \leq \KFC\zeta\para{\gamma}.
}
\end{remark}
\begin{definition}[$\omega$-smoothness]\label{def:omega}
Consider a function $\omega:\R_+\to\R_+$ such that $\omega\para{0}=0$ and $\xi\para{s}\eqdef\int_0^1\omega\para{st}dt$ is nondecreasing. A differentiable function $g:\HH\to\R$ is said to be $\omega$-smooth if, for every $x,y\in\HH$,
\newq{
\norm{\nabla g\para{x}- \nabla g\para{y}}{}\leq \omega\para{\norm{x-y}{}}
}
\end{definition}

\begin{remark}
	A classical consequence of $\omega$-smoothness is the following. If $g:\HH\to\R$ is $\omega$-smooth, for every $x,y\in\HH$ we have
\newq{
	f\para{y}\leq f\para{x} + \ip{\nabla f\para{x}, y-x}{} +  \xi \left(\norm{y-x}{}\right) \norm{y-x}{}.
}
\end{remark}

\begin{remark}
Note that being $\omega$-smooth is a stronger condition than being $\para{F,\zeta}$-smooth since every $\omega$-smooth function $f$ is also $\para{F,\zeta}$-smooth with $F = f$, $\zeta\para{t} = d_\C t\xi\para{d_\C t}$ and $\KFC\leq 1$. Additionally, the assumptions on $\xi$ being nondecreasing can be replaced by the sufficient condition that $\lim\limits_{t\to 0^+}\omega\para{t} = \omega\para{0}=0$.
\end{remark}
\ifdefined\COMPLETE
\else
\end{document}
\fi

\ifdefined\COMPLETE
\else
\documentclass[12pt]{article}
\input{package_header} 
\begin{document}
\fi


\section{Algorithm and Assumptions}\label{sec:alg}
For each $k\in\N$, we denote by $\Enk$ and $\snk$ random variables from $\prspace$ to $\HH_p$ and $\R_+$ respectively. In this context, $\Enk$ will represent the error in the gradient or proximal terms and $\snk$ will represent the error in the linear minimization oracle itself.\\
\begin{algorithm}[H]
    \SetAlgoLined
    \KwIn{$x_0\in \C \eqdef \dom\para{h}$; $\mu_0 \in \ran(A)$; $\seq{\gamma_k}$, $\seq{\beta_k}$, $\seq{\theta_k}, \seq{\rho_k} \in \ell_+$.}
    $k = 0$\\
    \Repeat{convergence}{
    $y_{k} =  \prox_{\beta_k g}\para{T x_k}$\\
    \vspace{0.25cm}
    $z_{k} = \nabla f(x_k)+T^*\left(Tx_k-y_k\right)/\beta_k+A^*\mu_k+\rho_kA^*\left(Ax_k-b\right) + \Enk$\\
    \vspace{0.2cm}
    $s_{k} \in \Argmin_{s\in \HH_p}\brac{h\para{s} + \ip{z_k,s}{}}$\\
    \vspace{0.2cm}
    $\nsk\in\brac{s\in \HH_p: h\para{s} + \ip{z_k,s}{} \leq h\para{s_k}+\ip{z_k,s_k}{} + \snk}$\\
    \vspace{0.2cm}
    $x_{k+1} = x_k -\gamma_k\para{x_k-\nsk}$\\
    \vspace{0.2cm}
    $\mu_{k+1} = \mu_k + \theta_k\para{Ax_{k+1}-b}$\\
    \vspace{0.2cm}
    $k \leftarrow k+1$\\
    \vspace{0.2cm}
    }
    \KwOut{$x_{k+1}$.}
\caption{Inexact Conditional Gradient with Augmented Lagrangian and Proximal-step (\icgalp)}
\label{alg:CGAL}
\end{algorithm}
To improve readability, we list some notation for the functionals we will employ throughout the analysis of the algorithm,
\nnewq{
\Phi \left(x\right) &\eqdef f\para{x} + g\para{Tx} + h\para{x};\\
\LL{x,\mu} &\eqdef f\para{x} + g\para{Tx} + h\para{x} + \ip{\mu,Ax-b}{};\\
\Lk{x,\mu} &\eqdef f\para{x} + g^{\beta_k}\para{Tx} + h\para{x} + \ip{\mu,Ax-b}{} + \frac{\rho_k}{2}\norm{Ax-b}{}^2;\\
\Ek{x,\mu} &\eqdef f\para{x} + g^{\beta_k}\para{Tx} + \ip{\mu,Ax-b}{} + \frac{\rho_k}{2}\norm{Ax-b}{}^2.
}
We can recognize $\LL{x,\mu}$ as the classical Lagrangian, $\Lk{x,\mu}$ as the augmented Lagrangian with smoothed $g$, and $\Ek{x,\mu}$ as the smooth part of $\Lk{x,\mu}$. With this notation in mind, we can see $z_k$ as $\nabla_x\Ek{x_k,\mu_k}$ and $\Enk$ as the error in the computation of $\grEk{x_k,\mu_k}$.

We define the filtration $\Filts\eqdef \seq{\filts_k}$ where $\filts_k \eqdef \sigma\para{x_0, \mu_0, \nsz,\ldots,\nsk}$ is the $\sigma$-algebra generated by the random variables $x_0, \mu_0, \nsz, \ldots, \nsk$. Furthermore, due to the error terms being contained in the direction finding step, we have that $x_{k+1}$ and $\mu_{k+1}$ are completely determined by $\filts_{k}$. Another noteworthy consequence of the error terms being contained in the direction finding step is that the primal iterates $\seq{x_k}$ remain in $\C$, as in the classical Frank-Wolfe algorithm, while the dual iterates $\seq{\mu_k}$ remain in $\ran\para{A}$.

\begin{subsection}{Assumptions}
\begin{subsubsection}{Assumptions on the functions}
We impose the following assumptions on the problem we consider; for some results, only a subset of them will be necessary:
\begin{enumerate}[label=(A.\arabic*)]
	\item \label{ass:A1} $f,g\circ T$, and $h$ belong to $\Gamma_0\para{\mc{H}_p}$
	\item \label{ass:f} The pair $\left(f,\C\right)$ is $\left(F,\zeta\right)$-smooth (see Definition \ref{def_smooth}), where we recall $\C\eqdef \dom\para{h}$
	\item \label{ass:compact} $\C$ is weakly compact (and thus contained in a ball of radius $\Ccon>0$)
	\item \label{ass:interior}$T\C \subset \dom(\p g)$ and $\sup\limits_{x\in \C}\norm{\sbrac{\p g\para{Tx}}^0}{} =\gcon< \infty$
	\item \label{ass:lip} $h$ is Lipschitz continuous relative to its domain $\C$ with constant $L_h \geq 0$, i.e., $\allowbreak \forall (x,z) \in \C^2$, $|h(x)-h(z)| \leq L_h \norm{x-z}{}$.
	\item \label{ass:existence} There exists a saddle-point $\para{\xs,\mus}\in\HH_p\times\HH_d$ for the Lagrangian $\mc{L}$
	\item \label{ass:rangeA} $\ran(A)$ is closed
	\item \label{ass:coercive} One of the following holds:
	\begin{enumerate}[label=(\alph*)]
	\item \label{ass:coerciveinfdim} $A^{-1}\para{b}\cap\inte\para{\dom\para{g\circ T}}\cap\inte\para{\C} \neq \emptyset$, where $A^{-1}\para{b}$ is the pre-image of $b$ under $A$
	\item \label{ass:coercivefindim} $\HH_p$ and $\HH_d$ are finite-dimensional and
	\nnewq{\label{eq:findimass}
	\begin{cases}
	A^{-1}\para{b} \cap \ri\para{\dom \para{g\circ T}} \cap \ri\para{\C}\neq\emptyset \\
	\qandq \\
	\ran\para{A^*}\cap \LinHull\para{\dom \para{g\circ T} \cap \C}^\bot = \brac{0} .
	\end{cases}
	}
	\end{enumerate}
\end{enumerate}

\end{subsubsection}

\begin{subsubsection}{Assumptions on the parameters and error terms}
We impose the following assumptions on the parameters and error terms and, as with the assumptions above, for some results only a subset will be necessary:
\begin{enumerate}[label=(P.\arabic*)]
	\item \label{ass:P1} $\seq{\gamma_k}\subset]0,1]$ and the sequences $\seq{\zeta\left(\gamma_k\right)}, \seq{\gamma_k^2/\beta_k}$ and $\seq{\gamma_k\beta_k}$ belong to $\ell^1_+$
	\item\label{notell1} $\seq{\gamma_k} \notin \ell^1$
	\item \label{ass:beta} $\seq{\beta_k} \in \ell_+$ is nonincreasing and converges to $0$
	\item \label{ass:rhodec}\label{ass:rhobound} $\seq{\rho_k} \in \ell_+$ is nondecreasing with $0 < \rhoinf \leq \rho_k \leq \rhosup<\pinfty$
	\item \label{ass:gkgkp1} For some positive constants $\gamconinf$ and $\gamcon$, $\gamconinf \leq \para{\gamma_k / \gamma_{k+1}} \leq \gamcon$
	\item \label{ass:thetagamma} $\seq{\theta_k}$ satisfies $\theta_k =  \frac{\gamma_k}{c}$ for some $c>0$ such that $\frac{\gamcon}{c}-\frac{\rhoinf}{2}< 0$
	\item \label{ass:cond} $\seq{\gamma_k}$ and $\seq{\rho_k}$ satisfy $\rho_{k+1}-\rho_k-\gamma_{k+1}\rho_{k+1}+\frac{2}{c}\gamma_k-\frac{\gamma_k^2}{c}\leq \gamma_{k+1}$ for $c$ in \ref{ass:thetagamma}
	\item \label{ass:error}$\seq{\gamma_{k+1}\EX{\norm{\Enkp}{}}{\filts_{k}}}\in\ell^1_+\para{\Filts}$ and $\seq{\gamma_{k+1}\EX{\snkp}{\filts_{k}}}\in\ell^1_+\para{\Filts}$.
\end{enumerate}
\begin{remark}
We will also denote the gradient of $\mc{E}_k$ with errors as
\newq{
\ngrEk{x,\mu} \eqdef \grEk{x,\mu} + \Enk.
}
It is possible to further decompose the error term $\Enk$, for instance, into $\fnk - T^*\gnk/\beta_k$ where $\fnk$ is the error in computing $\nabla f\para{x_k}$ and $\gnk$ is the error in evaluating $\prox_{\beta_k g}\para{Tx_k}$. In this case, the condition $\seq{\gamma_{k+1}\EX{\norm{\Enkp}{}}{\filts_{k}}}\in\ell^1_+\para{\Filts}$ in \ref{ass:error} is sufficiently satisfied by demanding that $\seq{\gamma_{k+1}\EX{\norm{\fnkp}{}}{\filts_{k}}}\in\ell^1_+\para{\Filts}$ and $\seq{\frac{\gamma_{k+1}}{\beta_{k+1}}\EX{\norm{\gnkp}{}}{\filts_{k}}}\in\ell^1_+\para{\Filts}$.
\end{remark}
\end{subsubsection}
\end{subsection}

\ifdefined\COMPLETE
\else
\end{document}
\fi

\ifdefined\COMPLETE
\else
\documentclass[12pt]{article}
\input{package_header} 
\begin{document}
\fi
\section{Main Results}\label{sec:mainres}
\ifdefined\COMPLETE
\else
\documentclass[12pt]{article}
\input{package_header} 
\begin{document}
\fi
\subsection{Preparatory Results}\label{sec:auxres}
\begin{lemma}\label{DescLemma}
	Suppose \ref{ass:A1}, \ref{ass:f} and \ref{ass:P1} hold. For each $k\in\N$, define the quantity
	\begin{equation}\label{Lk}
	L_k\eqdef \frac{\norrm{T}^2}{\beta_k}+\norrm{A}^2\rho_k.
	\end{equation} 
	Then, for each $k\in\N$, we have the following inequality,
\begin{equation*}
	\begin{split}
		\Ek{x_{k+1},\mu_k} & \leq \Ek{x_k,\mu_k} + \scal{\nabla_x \Ek{x_k,\mu_k}}{x_{k+1}-x_k} +  D_F\para{x_{k+1},x_{k}} \\
				   & \qquad + \frac{L_k}{2}\norrm{x_{k+1}-x_k}^2. 
	\end{split}
\end{equation*}
	\end{lemma}

\begin{proof}
See \cite[Lemma~4.5]{silveti}
	\end{proof}
	
\begin{lemma}\label{lemma:lowerbound}
	Suppose \ref{ass:A1} and \ref{ass:f} hold. Then, for each $k\in\N$ and for every $x\in\HH_p$,
	\begin{equation*}
		\begin{split}
			\Ek{x,\mu_k} & \geq \Ek{x_k,\mu_k} + \scal{\nabla_x \Ek{x_k,\mu_k}}{x-x_k}+\frac{\rho_k}{2}\norrm{A(x-x_k)}^2.
		\end{split}
	\end{equation*}
\end{lemma}
\begin{proof}
See \cite[Lemma~4.6]{silveti}.
\end{proof}

\begin{lemma}\label{convfeas}
Assume that \ref{ass:compact} and \ref{ass:rhodec} hold. Let $\seq{x_k}$ be the sequence of primal iterates generated by \algref{alg:CGAL} and $\Filts=\seq{\filts_k}$ as before. Then, for each $k\in\N$, we have the following estimate,
\newq{
\frac{\rho_k}{2}\norm{Ax_k-b}{}^2 - \frac{\rho_{k+1}}{2}\EX{\norm{Ax_{k+1}-b}{}^2}{\filts_{k-1}} \leq \rhosup d_\C\norm{A}{}\para{\norm{A}{}R + \norm{b}{}}\gamma_k \Pas.
}
\end{lemma}
\begin{proof}
For each $k\in\N$, by convexity of the function $\frac{\rho_{k+1}}{2}\norm{A\cdot-b}{}^2$ and the assumption \ref{ass:rhodec} that $\seq{\rho_k}$ is nondecreasing, we have,
\newq{
\frac{\rho_k}{2}\norm{Ax_k-b}{}^2 - \frac{\rho_{k+1}}{2}\norm{Ax_{k+1}-b}{}^2 &\leq \frac{\rho_{k+1}}{2}\norm{Ax_k-b}{}^2 - \frac{\rho_{k+1}}{2}\norm{Ax_{k+1}-b}{}^2\\
&\leq \ip{\nabla \para{\frac{\rho_{k+1}}{2}\norm{A\cdot-b}{}^2\para{x_k},x_k-x_{k+1}}}{}\\
&= \rho_{k+1}\ip{Ax_{k}-b,A\para{x_k-x_{k+1}}}{}.
}
Recall that, for each $k\in\N$, $x_{k+1}=x_k-\gamma_k\para{x_k-\nsk}$ and take the expectation to find,
\newq{
\frac{\rho_k}{2}\norm{Ax_k-b}{}^2 - \EX{\frac{\rho_{k+1}}{2}\norm{Ax_{k+1}-b}{}^2}{\filts_{k-1}} & \leq \rhosup \gamma_k\EX{\ip{Ax_k-b, A\para{x_k-\nsk}}{}}{\filts_{k-1}}\\
&\leq \rhosup\gamma_kd_\C\norm{A}{}\para{\norm{A}{}R+\norm{b}{}},
}
where we have used the Cauchy-Schwartz inequality and the boundedness of $\C$, assumed in \ref{ass:compact}, in the last inequality.
\end{proof}

\begin{remark}
The above result still holds if we replace both $\rho_k$ and $\rho_{k+1}$ by the constant $2$ and shift the index by $1$, i.e., for each $k\in\N$,
\newq{
\norm{Ax_{k+1}-b}{}^2 - \EX{\norm{Ax_{k+2}-b}{}^2}{\filts_{k}}\leq 2d_\C\norm{A}{}\para{\norm{A}{}R+\norm{b}{}}\gamma_{k+1}\Pas
}
\end{remark}

\begin{lemma}\label{convopt}
Suppose that \ref{ass:A1}-\ref{ass:existence} hold. Let $\seq{x_k}$ be the sequence of primal iterates generated by \algref{alg:CGAL} and $\mus$ a solution, which exists by \ref{ass:existence}, of the dual problem. Then, for each $k\in\N$, we have the following estimate,
\newq{
\LL{x_k,\mus} - \EX{\LL{x_{k+1},\mus}}{\filts_{k-1}} \leq \gamma_k d_\C\para{\gcon\norm{T}{} + D + L_h + \norm{\mus}{}\norm{A}{}}\Pas.
}
\end{lemma}
\begin{proof}
We recall the proof from \cite[Lemma~4.7]{silveti} with a slight modification to account for the inexactness of the algorithm. Define $u_{k}\eqdef\left[\partial g(Tx_{k})\right]^0$ and recall that, by \ref{ass:interior} and the fact that for all $k\in\N$, $x_k\in\C$, we have $\norrm{u_{k}} \leq \gcon$. By $\ref{ass:A1}$, the function $\Phi\para{x}\eqdef f\para{x} + g\para{Tx} + h\para{x}$ is convex. Then\fekn, 
	\newq{
	\LL{x_{k},\mus}-\LL{x_{k+1},\mus} &=\Phi(x_{k})-\Phi(x_{k+1})+\scal{\mus}{A\para{x_{k}-x_{k+1}}}\\
	&\leq \scal{u_{k}}{T(x_{k}-x_{k+1})}+\scal{\nabla f(x_{k})}{x_{k}-x_{k+1}}\\
	&\quad+L_h\norrm{x_{k}-x_{k+1}}+\norrm{\mus}\ \norrm{A}\ \norrm{x_{k}-x_{k+1}},
	}
where we used the subdifferential inequality \eqref{subdiff} on $g$ and $f$, the $L_h$-Lipschitz continuity of $h$ relative to $\C$ (see~\ref{ass:lip}), and the Cauchy-Schwartz inequality on the inner product. Since\fekn $x_{k+1}=x_{k}+\gamma_{k} \left(\nsk-x_{k}\right)$, we obtain\fekn
	\newq{
	\LL{x_{k},\mus}-\LL{x_{k+1},\mus}	
	&\leq \gamma_{k} \Big(\scal{u_{k}}{T(x_{k}-\nsk)}+\scal{\nabla f(x_{k})}{x_{k}-\nsk}+ L_h\norrm{x_{k}-\nsk}\\
	&\ \ +\norrm{\mus}\ \norrm{A}\ \norrm{x_{k}-\nsk}\Big)
	}
Now take the expectation with respect to the filtration $\filts_{k-1}$, such that $x_k$ is completely determined, to get\fekn
	\newq{
	\LL{x_{k},\mus}-\EX{\LL{x_{k+1},\mus}}{\filts_{k-1}}&\leq \gamma_{k} \Big(\EX{\scal{u_{k}}{T(x_{k}-\nsk)}}{\filts_{k-1}}+\EX{\scal{\nabla f(x_{k})}{x_{k}-\nsk}}{\filts_{k-1}}\\
	&\ \ + L_h\EX{\norrm{x_{k}-\nsk}}{\filts_{k-1}}+\norrm{\mus}\ \norrm{A}\ \EX{\norrm{x_{k}-\nsk}}{\filts_{k-1}}\Big)\\
	&\leq \gamma_{k} \diam\para{\gcon\norrm{T}+D+L_h+\norrm{\mus}\ \norrm{A}},
	}
	where we have used the Cauchy-Schwartz inequality, the boundedness of the set $\C$ by \ref{ass:compact}, the boundedness of $u_k$ by $\gcon$ by \ref{ass:interior}, and denoted by $D$ the constant $D \eqdef \sup_{x\in\C}\norrm{\nabla f(x)} < +\infty$. Note that $D$ exists and is bounded since $f$ is differentiable on an open set $\C_0$ containing $\C$ by \ref{ass:f} and \defref{def_smooth}.
\end{proof}


\ifdefined\COMPLETE
\else
\end{document}
\fi

\ifdefined\COMPLETE
\else
\documentclass[12pt]{article}
\input{package_header} 
\begin{document}
\fi

\subsection{Asymptotic feasibility}\label{sec:asfeas}
\begin{lemma}[Feasibility estimate]\label{feaslem}
Suppose that \ref{ass:A1} - \ref{ass:interior} and \ref{ass:existence} all hold. Consider the sequence of iterates $\seq{x_k}$ generated by Algorithm~\ref{alg:CGAL} with parameters satisfying \ref{ass:P1} and \ref{ass:beta}-\ref{ass:thetagamma}. For each $k\in\N$, define the two quantities, $\Delta_k^p$ and $\Delta_k^d$ in the following way, 
\newq{
\Delta^p_k \eqdef  \mc{L}_{k}\para{x_{k+1},\mu_k}-\minL_{k}\para{\mu_k},\quad \Delta^d_k \eqdef  \minL-\minL_{k}\para{\mu_k},
}
where we have denoted $\minL_{k}\para{\mu_k}\eqdef \min_{x}\mc{L}_k\para{x,\mu_k}$ and $\minL \eqdef  \mc{L}\para{\xs,\mus}$. Furthermore\fekn denote the sum $\Delta_k\eqdef \Delta_k^p + \Delta_k^d$.
We then have\fekn
\newq{
\EX{\Delta_{k+1}}{\filtfe_{k}} - \Delta_k	&\leq -\gamma_{k+1}\para{\frac{\gamconinf}{c}\norm{A\minx_{k+1}-b}{}^2 +\delta\norm{A\para{x_{k+1}-\minx_{k+1}}}{}^2}+\gamma_{k+1}^2\frac{L_{k+1}}{2}\diam^2\\
								&\quad +\KFC\zeta\para{\gamma_{k+1}} +\frac{\beta_{k}-\beta_{k+1}}{2}\gcon^2 +\para{\rho_{k+1}-\rho_k}\para{\norm{A}{}^2 \Ccon^2+\norm{b}{}^2} \\
								&\quad  + \gamma_{k+1}\EX{\snkp}{\filtfe_{k}}+\diam\gamma_{k+1}\EX{\norm{\Enkp}{}}{\filtfe_{k}}.
}
\end{lemma}
\begin{proof}
The proof here is adapted from the analogous result found in \cite[Theorem~4.1]{silveti}. As before, the quantity $\Delta^p_k \geq 0$ and can be seen as a primal gap at iteration $k$ while $\Delta^d_k$ may be negative but is uniformly bounded from below by our assumptions (see \cite[Theorem~4.1]{silveti}). We denote a minimizer of $\mc{L}_k\para{x,\mu_k}$ by $\minx_k \in \Argmin\limits_{x\in\mc{H}_p}\mc{L}_k\para{x,\mu_k}$, which exists and belongs to $\C$ by \ref{ass:A1}-\ref{ass:compact}. We have\fekn
\begin{equation*}
\begin{split}
\Delta_{k+1}-\Delta_k & = \mc{L}_{k+1}\para{x_{k+2},\mu_{k+1}}-\mc{L}_k\para{x_{k+1},\mu_{k+1}}+\theta_k\norm{Ax_{k+1}-b}{}^2\\
& \quad +2\left[\mc{L}_k\para{\tilde{x}_{k},\mu_{k}}-\mc{L}_{k+1}\para{\tilde{x}_{k+1},\mu_{k+1}}\right].
\end{split}
\end{equation*}
Recall that $\minx_k \in \Argmin\limits_{x\in\mc{H}_p}\mc{L}_k\para{x,\mu_k}$, that $g^{\beta_k}\leq g^{\beta_{k+1}}$ due to \ref{ass:beta} and \propref{moreau}\ref{moreauclaim4}, and that $\rho_k\leq\rho_{k+1}$ by \ref{ass:rhodec}. Then\fekn
\begin{equation*}
\begin{split}
\mc{L}_k\para{\tilde{x}_{k},\mu_{k}}-\mc{L}_{k+1}\para{\tilde{x}_{k+1},\mu_{k+1}}&\leq \mc{L}_k\para{\tilde{x}_{k+1},\mu_{k}}-\mc{L}_{k+1}\para{\tilde{x}_{k+1},\mu_{k+1}}\\
&=\left[g^{\beta_{k}}-g^{\beta_{k+1}}\right]\para{T \tilde{x}_{k+1}}+\frac{1}{2}\left[\rho_{k}-\rho_{k+1}\right]\norm{A\tilde{x}_{k+1}-b}{}^2\\
& \quad +\ip{\mu_{k}-\mu_{k+1},A\minx_{k+1}-b}{}\\
& \leq -\theta_{k} \ip{Ax_{k+1}-b,A\minx_{k+1}-b}{},
\end{split}
\end{equation*}
where we have used the fact that $\mu_{k+1} = \mu_k + \theta_k \para{Ax_{k+1}-b}$ coming from Algorithm~\ref{alg:CGAL}.
So we get\fekn
\begin{equation*}
\begin{split}
\Delta_{k+1}-\Delta_k & \leq \mc{L}_{k+1}\para{x_{k+2},\mu_{k+1}}-\mc{L}_k\para{x_{k+1},\mu_{k+1}}+\theta_k\norm{Ax_{k+1}-b}{}^2\\
& \quad -2\theta_{k} \ip{Ax_{k+1}-b,A\minx_{k+1}-b}{}.
\end{split}
\end{equation*}
Note that\fekn
\newq{
	\mc{L}_k\para{x_{k+1},\mu_{k+1}}=\mc{L}_{k+1}\para{x_{k+1},\mu_{k+1}}-\left[g^{\beta_{k+1}}-g^{\beta_k}\right]\para{T x_{k+1}}-\para{\frac{\rho_{k+1}-\rho_k}{2}}\norm{Ax_{k+1}-b}{}^2.
}
Then\fekn
\newqml{
	\Delta_{k+1}-\Delta_k \leq \mc{L}_{k+1}\para{x_{k+2},\mu_{k+1}}-\mc{L}_{k+1}\para{x_{k+1},\mu_{k+1}}+g^{\beta_{k+1}}\para{T x_{k+1}}-g^{\beta_k}\para{T x_{k+1}}\\
	 +\para{\frac{\rho_{k+1}-\rho_k}{2}}\norm{Ax_{k+1}-b}{}^2+\theta_k\norm{Ax_{k+1}-b}{}^2-2\theta_k\ip{Ax_{k+1}-b,A\minx_{k+1}-b}{}.
}
We denote by $\Te \eqdef\mc{L}_{k+1}\para{x_{k+2},\mu_{k+1}}-\mc{L}_{k+1}\para{x_{k+1},\mu_{k+1}}$ and the remaining part of the right-hand side by $\Tee$.	For the moment, we focus our attention on $\Te$. Recall that $\mc{L}_{k}\para{x,\mu_{k}}=\mc{E}_k\para{x,\mu_{k}}+h\para{x}$ and apply \lemref{DescLemma} between points $x_{k+2}$ and $x_{k+1}$, to get\fekn
\newqml{
\Te \leq h\para{x_{k+2}}-h\para{x_{k+1}}+\ip{\nabla_x \mc{E}_{k+1}\para{x_{k+1},\mu_{k+1}},x_{k+2}-x_{k+1}}{}\\
+\frac{L_{k+1}}{2}\norm{x_{k+2}-x_{k+1}}{}^2+D_F\para{x_{k+2}, x_{k+1}}.
}
By \ref{ass:A1} we have that $h$ is convex and thus, since $x_{k+2}$ is a convex combination of $x_{k+1}$ and $\nskp$, we get\fekn
\newq{
\Te &\leq \gamma_{k+1}\para{h\para{\nskp}-h\para{x_{k+1}}+\ip{\nabla_x\mc{E}_{k+1}\para{x_{k+1},\mu_{k+1}}, \nskp-x_{k+1}}{}}\\
	&\quad+\frac{L_{k+1}}{2}\norm{x_{k+2}-x_{k+1}}{}^2+D_F\para{x_{k+2}, x_{k+1}}\\
	& =  \gamma_{k+1}\Big(h\para{\nskp}-h\para{x_{k+1}}+\ip{\ngrEkp{x_{k+1},\mu_{k+1}}, \nskp-x_{k+1}}{}\\
	&\quad+ \ip{\nabla_x\mc{E}_{k+1}\para{x_{k+1},\mu_{k+1}}-\ngrEkp{x_{k+1},\mu_{k+1}}, \nskp-x_{k+1}}{}\Big)\\
	&\quad+\frac{L_{k+1}}{2}\norm{x_{k+2}-x_{k+1}}{}^2+D_F\para{x_{k+2}, x_{k+1}}\\
	&= \gamma_{k+1}\Big(h\para{\nskp}-h\para{x_{k+1}}+\ip{\ngrEkp{x_{k+1},\mu_{k+1}}, \nskp-x_{k+1}}{}\\
	&\quad - \ip{\Enkp, \nskp-x_{k+1}}{}\Big)+\frac{L_{k+1}}{2}\norm{x_{k+2}-x_{k+1}}{}^2+D_F\para{x_{k+2}, x_{k+1}}
}
Applying the definition of $\nsk$ as the approximate minimizer of the linear minimization oracle gives\fekn
\newq{
\Te	&\leq \gamma_{k+1}\Big(h\para{s_{k+1}}-h\para{x_{k+1}}+\ip{\ngrEkp{x_{k+1},\mu_{k+1}}, s_{k+1}-x_{k+1}}{} + \snkp\\
	& \quad-\ip{\Enkp, \nskp-x_{k+1}}{}\Big)+\frac{L_{k+1}}{2}\norm{x_{k+2}-x_{k+1}}{}^2+D_F\para{x_{k+2}, x_{k+1}}.
}
Now we can apply the definition of $s_{k+1}$ as the minimizer of the linear minimization oracle and \lemref{lemma:lowerbound} to get\fekn
\newq{
\Te & \leq \gamma_{k+1}\Big(h\para{\minx_{k+1}}-h\para{x_{k+1}}+\ip{\ngrEkp{x_{k+1},\mu_{k+1}}, \minx_{k+1}-x_{k+1}}{} + \snkp\\
	& \quad-\ip{\Enkp, \nskp-x_{k+1}}{}\Big)+\frac{L_{k+1}}{2}\norm{x_{k+2}-x_{k+1}}{}^2+D_F\para{x_{k+2}, x_{k+1}}\\
	& = \gamma_{k+1}\Big(h\para{\minx_{k+1}}-h\para{x_{k+1}}+\ip{\nabla_x\mc{E}_{k+1}\para{x_{k+1},\mu_{k+1}}, x_{k+1}-\minx_{k+1}}{} + \snkp\\
	& \quad-\ip{\Enkp, \nskp-\minx_{k+1}}{}\Big)+\frac{L_{k+1}}{2}\norm{x_{k+2}-x_{k+1}}{}^2+D_F\para{x_{k+2}, x_{k+1}}\\
	&\leq \gamma_{k+1}\Big(h\para{\minx_{k+1}}-h\para{x_{k+1}}+\mc{E}_{k+1}\para{\minx_{k+1},\mu_{k+1}}-\mc{E}_{k+1}\para{x_{k+1},\mu_{k+1}}-\frac{\rho_{k+1}}{2}\norm{A\para{x_{k+1}-\minx_{k+1}}}{}^2\\
	& \quad + \snkp-\ip{\Enkp, \nskp-\minx_{k+1}}{}\Big)+\frac{L_{k+1}}{2}\norm{x_{k+2}-x_{k+1}}{}^2+D_F\para{x_{k+2}, x_{k+1}}\\
	&= \gamma_{k+1}\Big(\mc{L}_{k+1}\para{\minx_{k+1},\mu_{k+1}}-\mc{L}_{k+1}\para{x_{k+1},\mu_{k+1}}-\frac{\rho_{k+1}}{2}\norm{A\para{x_{k+1}-\minx_{k+1}}}{}^2 + \snkp\\
	&\quad -\ip{\Enkp, \nskp-\minx_{k+1}}{}\Big)+\frac{L_{k+1}}{2}\norm{x_{k+2}-x_{k+1}}{}^2+D_F\para{x_{k+2}, x_{k+1}} \\
	&\leq -\frac{\gamma_{k+1}\rho_{k+1}}{2}\norm{A\para{x_{k+1}-\minx_{k+1}}}{}^2 + \gamma_{k+1}\Big( \snkp +\ip{\Enkp,\minx_{k+1}-\nskp}{} \Big)\\
	&\quad +\frac{L_{k+1}}{2}\norm{x_{k+2}-x_{k+1}}{}^2+D_F\para{x_{k+2}, x_{k+1}} ,
}
where we used that $\minx_{k+1}$ is a minimizer of $\mc{L}_{k+1}\para{\cdot,\mu_{k+1}}$ in the last inequality. Now combining $\Te$ and $\Tee$ and using the Pythagoras identity we have\fekn
\newq{
\Delta_{k+1}-\Delta_k	&\leq -\theta_{k}\norm{A\minx_{k+1}-b}{}^2+\para{\theta_{k}-\gamma_{k+1}\frac{\rho_{k+1}}{2}}\norm{A\para{x_{k+1}-\minx_{k+1}}}{}^2\\
						&\quad +\frac{L_{k+1}}{2}\norm{x_{k+2}-x_{k+1}}{}^2+D_F\para{x_{k+2}, x_{k+1}}+\left[g^{\beta_{k+1}}-g^{\beta_k}\right]\para{T x_{k+1}} \\
						&\quad +\frac{\rho_{k+1}-\rho_k}{2}\norm{Ax_{k+1}-b}{}^2+ \gamma_{k+1}\Big( \snkp +\ip{\Enkp,\minx_{k+1}-\nskp}{} \Big).
}
\ \\
Now take the expectation with respect to $\filtfe_k = \filts_k = \sigma\para{x_0, \mu_0,\nsz,\ldots,\nsk}$, which completely determines $x_{k+1}$, $\minx_{k+1}$, and $\mu_{k+1}$. We are also going to perform the following estimations.
\begin{itemize}
	\item Under \ref{ass:gkgkp1} and \ref{ass:thetagamma}, we have that\fekn $\theta_k=\gamma_{k}/c$ with $\gamconinf \gamma_{k+1} \leq \gamma_k$ and so that
	$$-\theta_{k} \leq -\tfrac{\gamconinf }{c}\gamma_{k+1}. $$
	\item Again by \ref{ass:thetagamma}, we have\fekn $\theta_k=\gamma_{k}/c$ for some $c>0$ such that
	\newq{
		\exists \delta>0,\quad \frac{\gamcon}{c} -\frac{\rhoinf}{2}=-\delta < 0,
	}
	where $\gamcon$ is the constant such that\fekn $\gamma_{k}\leq \gamcon \gamma_{k+1}$ (see \ref{ass:gkgkp1}).
	Then, using again \ref{ass:gkgkp1} and the above inequality\fekn
	\nnewq{\label{eq:aux1}
		\theta_{k}-\gamma_{k+1}\frac{\rho_{k+1}}{2}\leq \para{\frac{\gamcon}{c}-\frac{\rho_{k+1}}{2}}\gamma_{k+1}\leq\para{\frac{\gamcon}{c}-\frac{\rhoinf}{2}}\gamma_{k+1}=-\delta\gamma_{k+1}.
	}
\item By Algorithm~\ref{alg:CGAL}\fekn $x_{k+2} - x_{k+1} = \gamma_{k+1}\para{\nskp-x_{k+1}}$. Since $\nskp$ and $x_{k+1}$ are both in $\C$ and $\C$ is bounded due to \ref{ass:compact}\fekn
$$\frac{L_{k+1}}{2}\EX{\norm{x_{k+2}-x_{k+1}}{}^2}{\filtfe_{k}} = \frac{L_{k+1}}{2}\gamma^2_{k+1}\EX{\norm{\nskp-x_{k+1}}{}^2}{\filtfe_{k}} \leq \frac{L_{k+1}}{2}\gamma^2_{k+1}\diam^2.$$
\item Recall that, by \ref{ass:f}, $f$ is $\para{F,\zeta}$-smooth and invoke \remref{fremark}, to get
$$\EX{D_F\para{x_{k+2}, x_{k+1}}}{\filtfe_{k}}\leq \KFC\zeta\para{\gamma_{k+1}}.$$
\item By \propref{moreau}\ref{moreauclaim4}  and assumption \ref{ass:interior},
$$\EX{\left[g^{\beta_{k+1}}-g^{\beta_k}\right]\para{T x_{k+1}}}{\filtfe_{k}}\leq \frac{\beta_k-\beta_{k+1}}{2}\EX{\norm{\sbrac{\partial g\para{Tx_{k+1}}}^0}{}^2}{\filtfe_{k}}\leq \frac{\beta_k-\beta_{k+1}}{2}M^2.$$
\item We also have, using Jensen's inequality and \ref{ass:compact}\fekn $$\para{\frac{\rho_{k+1}-\rho_k}{2}}\EX{\norm{Ax_{k+1}-b}{}^2}{\filtfe_{k}}\leq \para{\rho_{k+1}-\rho_k}\para{\norm{A}{}^2 \Ccon^2+\norm{b}{}^2}.$$
\end{itemize}

%
%
%
%
\ \\
In total\fekn
\newq{
\EX{\Delta_{k+1}}{\filtfe_{k}} - \Delta_k	&\leq -\tfrac{\gamconinf }{c}\gamma_{k+1}\norm{A\minx_{k+1}-b}{}^2-\delta\gamma_{k+1}\norm{A\para{x_{k+1}-\minx_{k+1}}}{}^2\\
								&\quad +\frac{L_{k+1}}{2}\gamma^2_{k+1}\diam^2+\KFC\zeta\para{\gamma_{k+1}} \\
								&\quad +\frac{\beta_k-\beta_{k+1}}{2}M^2 +\para{\rho_{k+1}-\rho_k}\para{\norm{A}{}^2 \Ccon^2+\norm{b}{}^2}\\
								&\quad +\gamma_{k+1}\Big( \EX{\snkp}{\filtfe_{k}} +\EX{\ip{\Enkp,\minx_{k+1}-\nskp}{}}{\filtfe_{k}} \Big).
}
Using Cauchy-Schwarz together with the fact that $\minx_{k+1}$ and $\nskp$ are in $\C$, which is bounded by \ref{ass:compact}, we also have\fekn
\nnewq{
\gamma_{k+1}\EX{\ip{\Enkp,\minx_{k+1}-\nskp}{}}{\filtfe_{k}} \leq
\gamma_{k+1}\diam\EX{\norm{\Enkp}{}}{\filtfe_k},
}
which gives\fekn
\nnewq{\label{eq:estfeas}
\EX{\Delta_{k+1}}{\filtfe_{k}} - \Delta_k	&\leq -\frac{\gamconinf}{c}\gamma_{k+1}\norm{A\minx_{k+1}-b}{}^2-\delta\gamma_{k+1}\norm{A\para{x_{k+1}-\minx_{k+1}}}{}^2+\gamma_{k+1}^2\frac{L_{k+1}}{2}\diam^2\\
								&\quad +\KFC\zeta\para{\gamma_{k+1}} +\frac{\beta_{k}-\beta_{k+1}}{2}\gcon^2 +\para{\rho_{k+1}-\rho_k}\para{\norm{A}{}^2 \Ccon^2+\norm{b}{}^2}\\
								&\quad   + \gamma_{k+1}\EX{\snkp}{\filtfe_{k}} +\gamma_{k+1}\diam\EX{\norm{\Enkp}{}}{\filtfe_k},
}
and we conclude by trivial manipulations.
\end{proof}

\begin{theorem}[Feasibility]\label{feastheorem}
Suppose that \ref{ass:A1}-\ref{ass:interior} and \ref{ass:existence} all hold. For a sequence $\seq{x_k}$ generated by Algorithm~\ref{alg:CGAL} using parameters satisfying \ref{ass:P1} - \ref{ass:thetagamma} and \ref{ass:error} we have,
\begin{enumerate}[label=(\roman*)]
	\item \label{feaslim}Asymptotic feasbility: $\lim\limits_{k\to\infty} \norm{Ax_k-b}{} =0$ $\Pas$
	\item \label{feaspwrate}Pointwise rate: \begin{equation}\label{pointratefeas}
	\begin{split}
	\inf_{0 \leq i \leq k} \norrm{Ax_{i}-b} = O\para{\frac{1}{\sqrt{\Gamma_{k}}}} \Pas \qandq\\ \text{$\exists$ a subsequence $\subseq{x_{k_j}}$ such that } \norrm{Ax_{k_j}-b} \leq \frac{1}{\sqrt{\Gamma_{k_j}}} \Pas.
	\end{split}
	\end{equation}
where $\Gamma_k \eqdef \sum_{i=0}^{k} \gamma_i$.
	\item \label{feasergrate}Ergodic rate: let $\avx_k \eqdef \sum_{i=0}^{k} \gamma_i x_i/\Gamma_k$. Then
\begin{equation}\label{ergratefeas}
	\begin{split}
	\norrm{A\avx_k-b} = O\para{\frac{1}{\sqrt{\Gamma_{k}}}} \Pas.
	\end{split}
	\end{equation}
\end{enumerate}
\end{theorem}
\begin{proof}
Our goal is to first apply \lemref{combstoch} and then apply \lemref{barty}. By \lemref{feaslem}, we have\fekn
\nnewq{\label{feasest}
\EX{\Delta_{k+1}}{\filtfe_{k}} - \Delta_k	&\leq -\gamma_{k+1}\para{\frac{\gamconinf}{c}\norm{A\minx_{k+1}-b}{}^2 +\delta\norm{A\para{x_{k+1}-\minx_{k+1}}}{}^2}+\gamma_{k+1}^2\frac{L_{k+1}}{2}\diam^2\\
								&\quad +\KFC\zeta\para{\gamma_{k+1}} +\frac{\beta_{k}-\beta_{k+1}}{2}\gcon^2 +\para{\rho_{k+1}-\rho_k}\para{\norm{A}{}^2 \Ccon^2+\norm{b}{}^2} \\
								&\quad  + \gamma_{k+1}\EX{\snkp}{\filtfe_{k}}+\diam\gamma_{k+1}\EX{\norm{\Enkp}{}}{\filtfe_{k}}.
}
Because of \ref{ass:P1} and \ref{ass:rhobound}, and in view of the definition of $L_{k+1}$ in \eqref{Lk}, we have the following,
\newq{
\seq{\frac{L_{k+1}}{2}\gamma_{k+1}^2\diam^2}&= \seq{\frac{1}{2}\para{\frac{\norm{T}{}^2}{\beta_{k+1}}+\norm{A}{}^2\rho_{k+1}}\gamma_{k+1}^2\diam^2} \in\ell^1_+.
}
For the telescopic terms from the right hand side of \eqref{feasest} we have
\newq{
\seq{\frac{\beta_{k}-\beta_{k+1}}{2}\gcon^2} \in\ell^1_+ \  \tandt \ \seq{\para{\rho_{k+1}-\rho_k}\para{\norm{A}{}^2 \Ccon^2+\norm{b}{}^2}}\in\ell^1_+,}

where $\Ccon$ is the constant arising from \ref{ass:compact}.
Under \ref{ass:P1} we also have that
\newq{
\seq{K_{\para{F,\zeta,\C}}\zeta\para{\gamma_{k+1}}}\in\ell^1_+.
}
Finally, due to \ref{ass:error}, we also have
\newq{
\seq{\gamma_{k+1}\EX{\snkp}{\filtfe_k}}\in\ell^1_+\para{\Filtfe},\quad \seq{\diam\gamma_{k+1}\EX{\norm{\Enkp}{}}{\filtfe_{k}}}\in\ell^1_+\para{\Filtfe}.
}
Using the notation of \lemref{combstoch}, we set\fekn
\newq{
r_k &= \Delta_k, \quad a_k = \gamma_{k+1}\para{\frac{\gamconinf}{c}\norm{A\minx_{k+1} -b}{}^2 + \delta\norm{A\para{x_{k+1}-\minx_{k+1}}}{}^2}, \tandt\\
z_k &= \frac{L_{k+1}}{2}\gamma_{k+1}^2\diam^2+\KFC\zeta\para{\gamma_{k+1}}+\frac{\beta_{k}-\beta_{k+1}}{2}\gcon^2+\para{\frac{\rho_{k+1}-\rho_k}{2}}\norm{Ax_{k+1}-b}{}^2\\
	&\quad +\gamma_{k+1}\EX{\snkp}{\filtfe_k}+\diam\gamma_{k+1}\EX{\norm{\Enkp}{}}{\filtfe_{k}}.
}
We have shown above that \fekn
\[
\EX{r_{k+1}}{\filtfe_k} - r_{k}\leq  - a_k + z_k ,
\]
where $\seq{z_k} \in \ell_+^1\para{\Filtfe}$, and $r_k$ is bounded from below. We then deduce using \lemref{combstoch} that $\seq{r_k}$ is convergent $\Pas$ and 
\nnewq{\label{eq:starl1}
\seq{\gamma_{k}\norm{A\minx_{k}-b}{}^2} \in \ell^1_+\para{\Filtfe}, \quad \seq{\gamma_{k}\norm{A\para{x_{k}-\minx_{k}}}{}^2}\in\ell^1_+\para{\Filtfe}.
}
Consequently, 
\nnewq{\label{eq:feasxkl1}
\seq{\gamma_{k}\norm{Ax_{k}-b}{}^2}\in\ell^1_+\para{\Filtfe},
}
since by the Cauchy-Schwarz inequality,
\newq{
\sum\limits_{k=1}^\infty \gamma_{k}\norm{Ax_{k}-b}{}^2 &\leq 2\sum\limits_{k=1}^\infty \gamma_{k}\para{\norm{A\para{x_{k}-\minx_{k}}}{}^2+\norm{A\minx_{k}-b}{}^2} < \pinfty.
}
To finish proving \ref{feaslim} we simply apply \lemref{convfeas} (with the remark which follows) and the conditions of \lemref{barty} are satisfied. Then, \ref{feaspwrate} and \ref{feasergrate} follow directly from the results of \cite[Theorem~4.1]{silveti}.
\end{proof}

\ifdefined\COMPLETE
\else
\end{document}
\fi

\ifdefined\COMPLETE
\else
\documentclass[12pt]{article}
\input{package_header} 
\begin{document}
\fi

\subsection{Optimality}\label{sec:opt}
The following lemmas regard the boundedness of the sequence of dual iterates $\seq{\mu_k}$ and the uniform boundedness of the Lagrangian. They were shown in the deterministic setting in \cite{silveti} and trivially extend to the stochastic case in light of \thmref{feastheorem}.
\begin{lemma}
Suppose that \ref{ass:A1}-\ref{ass:compact}, \ref{ass:existence}-\ref{ass:coercive}, and \ref{ass:P1}-\ref{ass:thetagamma} all hold. Then the sequence of dual iterates $\seq{\mu_k}$ generated by \algref{alg:CGAL} is bounded.
\end{lemma}
\begin{proof}
See \cite[Lemma~4.9]{silveti}.
\end{proof}

\begin{lemma}\label{unifbound}
Under \ref{ass:A1}-\ref{ass:coercive} and \ref{ass:P1}-\ref{ass:thetagamma}, the composite function $f+g\circ T + h$ is uniformly bounded on $\C$ and we have
\nnewq{
\Mcstopt\eqdef \sup\limits_{x\in\C}\absv{f\para{x}+g\para{Tx}+h\para{x}} + \sup\limits_{k\in\N}\norm{\mu_k}{}\para{\norm{A}{}R+b}<+\infty,
}
where $R$ is the radius from \ref{ass:compact}.
\end{lemma}
\begin{proof}
The proof follows directly from \cite[Lemma~4.10]{silveti} with the addition of \thmref{feastheorem}.
\end{proof}

We now begin with the main energy estimate needed to show the convergence of the Lagrangian values to optimality.
\begin{lemma}[Optimality estimate]\label{optlem}
Recall the constants $c$, $L_k$, $\gcon$, $D$, and $L_h$ from \ref{ass:thetagamma}, \lemref{DescLemma}, \ref{ass:interior}, \lemref{convfeas}, and \ref{ass:lip}, respectively.
Define\fekn
\newq{\label{r}
r_k \eqdef \para{1-\gamma_k}\Lk{x_x,\mu_k} + \frac{c}{2}\norm{\mu_k-\mus}{}^2
}
and
\newq{
C_k \eqdef \frac{L_k}{2}\diam^2 + \diam\para{\gcon\norrm{T}+D+L_h+\norrm{\mus}\ \norrm{A}}.
} 
Then, under \ref{ass:A1}-\ref{ass:coercive} and \ref{ass:P1}-\ref{ass:cond}, for the sequences $\seq{x_k}$ and $\seq{\mu_k}$ generated by Algorithm~\ref{alg:CGAL}, using the filtration $\Filtop = \seq{\filtop_k}$ with $\filtop_k = \filts_{k-1}$, the following inequality holds\fekn
\nnewq{
	\begin{split}
\EX{r_{k+1}}{\filtop_{k}} - r_k  &\leq -\gamma_k \para{\LL{x_{k},\mus}-\LL{\xs,\mus} + \frac{\rho_k}{2}\norm{Ax_k-b}{}^2} + \frac{\gamma_{k+1}}{2}\EX{\norm{Ax_{k+1}-b}{}^2}{\filtop_{k}}\\
	& \ \  +\para{\beta_k-\beta_{k+1}}\frac{\gcon^2}{2}+ \para{\gamma_k-\gamma_{k+1}}\Mcstopt+ \gamma_k\beta_{k}\frac{\gcon^2}{2}+ \KFC\zeta\para{\gamma_k}+\gamma_k^2C_k\\
	& \ \ \ +\diam\gamma_k\EX{\norm{\Enk}{}}{\filtop_{k}} +\gamma_k\EX{\snk}{\filtop_{k}} \Pas.
	\end{split}
}
\end{lemma}
\begin{proof}
Applying \lemref{lemma:lowerbound} to the points $\xs$ and $x_k$ we have\fekn
\newq{
\Ek{\xs,\mu_k} 	& \geq \Ek{x_k,\mu_k} + \scal{\grEk{x_k,\mu_k}}{\xs-x_k}+\frac{\rho_k}{2}\norm{A(\xs-x_k)}{}^2\\
				& = \Ek{x_k,\mu_k} + \ip{\ngrEk{x_k,\mu_k}, \xs-x_k}{} + \ip{\Enk, x_k-\xs}{} +\frac{\rho_k}{2}\norm{A(\xs-x_k)}{}^2\\
				& = \Ek{x_k,\mu_k} + \ip{\ngrEk{x_k,\mu_k}, \xs-x_k}{} + h\para{\xs} - h\para{\xs} + \ip{\Enk, x_k-\xs}{}\\
				&\quad +\frac{\rho_k}{2}\norm{A(\xs-x_k)}{}^2.
}
By the definition of $s_k$ as a minimizer and the definition of $\nsk$ we further have\fekn
\nnewq{\label{eq:strconv}
\Ek{\xs,\mu_k}	& \geq \Ek{x_k,\mu_k} + \ip{\ngrEk{x_k,\mu_k}, s_k-x_k}{} + h\para{s_k} - h\para{\xs} + \ip{\Enk, x_k-\xs}{}\\
				&\quad +\frac{\rho_k}{2}\norm{A(\xs-x_k)}{}^2\\
				& \geq \Ek{x_k,\mu_k} + \ip{\ngrEk{x_k,\mu_k}, \nsk-x_k}{} + h\para{\nsk} - \snk- h\para{\xs} + \ip{\Enk, x_k-\xs}{}\\
				&\quad +\frac{\rho_k}{2}\norm{A(\xs-x_k)}{}^2.
}
%
%
%
%
%
%
From Lemma \ref{DescLemma} applied to the points $x_{k+1}$ and $x_k$ and by definition of $x_{k+1}\eqdef x_k+\gamma_k\para{\nsk-x_k}$ in Algorithm~\ref{alg:CGAL}, we also have\fekn
\begin{equation*}
	\begin{split}
		\Ek{x_{k+1},\mu_k} & \leq \Ek{x_k,\mu_k} + \scal{\grEk{x_k,\mu_k}}{x_{k+1}-x_k} +  D_F\para{x_{k+1},x_k} + \frac{L_k}{2}\norm{x_{k+1}-x_k}{}^2\\
		&= \Ek{x_k,\mu_k} + \gamma_k \scal{\grEk{x_k,\mu_k}}{\nsk-x_k} +D_F\para{x_{k+1},x_k}+\gamma_k^2\frac{L_k}{2}\norm{\nsk-x_k}{}^2\\
		&= \Ek{x_k,\mu_k} + \gamma_k \scal{\ngrEk{x_k,\mu_k}}{\nsk-x_k} +\gamma_k \ip{\Enk, x_k-\nsk}{}+D_F\para{x_{k+1},x_k}\\
		&\quad +\gamma_k^2\frac{L_k}{2}\norm{\nsk-x_k}{}^2.
	\end{split}
\end{equation*}
We combine the latter with \eqref{eq:strconv}, to get\fekn
\nnewq{\label{eq:strconv+desc}
\Ek{x_{k+1},\mu_k}				&\leq \Ek{x_k,\mu_k}+\gamma_k \ip{\Enk, \xs -\nsk}{} +D_F\para{x_{k+1},x_k}+\gamma_k^2\frac{L_k}{2}\norm{\nsk-x_k}{}^2\\
								&\quad +\gamma_k \para{\Ek{\xs,\mu_k}+h(\xs)-\Ek{x_k,\mu_k}-h(\nsk)-\frac{\rho_k}{2}\norm{Ax_k-b}{}^2 + \snk}.
}
By convexity of $h$ from \ref{ass:A1} and the definition of $x_{k+1}$, we have\fekn
\nnewq{\label{eq:hconv}
\Lk{x_{k+1},\mu_k} - \Lk{x_k,\mu_k}	&= \Ek{x_{k+1},\mu_k} - \Ek{x_k,\mu_k} + h\para{x_{k+1}}-h\para{x_k}\\
									&\leq \Ek{x_{k+1},\mu_k} - \Ek{x_k,\mu_k} + \gamma_k\para{h\para{\nsk} - h\para{x_k}} 
}
Combining \eqref{eq:strconv+desc} and \eqref{eq:hconv}, we obtain\fekn
\nnewq{\label{eq:strconv+desc+hconv}
\Lk{x_{k+1},\mu_k} - \Lk{x_k,\mu_k}	&\leq \gamma_k\para{\Ek{\xs,\mu_k}+h(\xs)-\Ek{x_k,\mu_k} - h\para{x_k}}+D_F\para{x_{k+1},x_k}+\\
								&\quad \gamma_k^2\frac{L_k}{2}\norm{\nsk-x_k}{}^2+\gamma_k\para{ \ip{\Enk, \xs -\nsk}{} -\frac{\rho_k}{2}\norm{Ax_k-b}{}^2 +\snk}\\
								&= \gamma_k\para{\Lk{\xs,\mu_k} - \Lk{x_k,\mu_k}} + D_F\para{x_{k+1},x_k} + \gamma_k^2\frac{L_k}{2}\norm{\nsk-x_k}{}^2\\
								&\quad +\gamma_k\para{ \ip{\Enk, \xs -\nsk}{} -\frac{\rho_k}{2}\norm{Ax_k-b}{}^2 +\snk}
}
Recalling the definition of $\mu_{k+1}\eqdef \mu_k + A\para{x_{k+1}-b}$ in Algorithm~\ref{alg:CGAL}, we have\fekn
\newq{
	\Lk{x_{k+1},\mu_{k+1}}-\Lk{x_{k+1},\mu_{k}}=\scal{\mu_{k+1}-\mu_k}{Ax_{k+1}}=\theta_k\norm{Ax_{k+1}-b}{}^2.
}
We combine the above and \eqref{eq:strconv+desc+hconv} to get\fekn
\nnewq{\label{eq:estt}
\Lk{x_{k+1},\mu_{k+1}} -\Lk{x_k,\mu_k}
					&\leq\theta_k\norm{Ax_{k+1}-b}{}^2 +\gamma_k\para{\Lk{\xs,\mu_k} - \Lk{x_k,\mu_k}} + D_F\para{x_{k+1},x_k}\\
					&\quad + \gamma_k^2\frac{L_k}{2}\norm{\nsk-x_k}{}^2 +\gamma_k\para{ \ip{\Enk, \xs -\nsk}{} -\frac{\rho_k}{2}\norm{Ax_k-b}{}^2 +\snk}.
}
Notice that the update of the dual variable $\mu$ can be interpreted as a prox operator in the following way,
\newq{
\mu_{k+1} = \argmin\limits_{\mu\in \HH_d}\brac{-\Lk{x_{k+1},\mu} + \frac{1}{2\theta_k}\norm{\mu-\mu_k}{}^2}.
}
Then, using \lemref{goodineq}, we get\fekn
\nnewq{\label{eq:proxmu}
0	&\geq \theta_k\para{\Lk{x_{k+1},\mus} - \Lk{x_{k+1},\mu_{k+1}}} + \frac{1}{2}\para{\norm{\mu_{k+1}-\mus}{}^2 - \norm{\mu_k-\mus}{}^2 + \norm{\mu_{k+1}-\mu_k}{}^2}\\
	&= \theta_k\para{\Lk{x_{k+1},\mus} - \Lk{x_{k+1},\mu_{k+1}}} + \frac{1}{2}\para{\norm{\mu_{k+1}-\mus}{}^2 - \norm{\mu_k-\mus}{}^2 + \theta_k^2\norm{Ax_{k+1}-b}{}^2}.
}
	Recall that, by \ref{ass:thetagamma}, $\theta_k =  \gamma_k/c$. Multiply \eqref{eq:proxmu} by $c$ and sum with \eqref{eq:estt}, to obtain\fekn
	\begin{eqnarray*}
		& (1-c\theta_k)\Lk{x_{k+1},\mu_{k+1}} -(1-c\theta_k)\Lk{x_k,\mu_k} + \frac{c}{2}\para{\norm{\mu_{k+1}-\mus}{}^2-\norm{\mu_k-\mus}{}^2} \\
		& \leq \ \ \para{\theta_k -\frac{c\theta_k^2}{2}}\norm{Ax_{k+1}-b}{}^2+\gamma_k \para{\Lk{\xs,\mu_k}-\Lk{x_k,\mu_k}}-c\theta_k \para{\Lk{x_{k+1},\mu}-\Lk{x_k,\mu_k}} \\
		& \ \ \ -\frac{\rho_k\gamma_k}{2}\norm{Ax_{k}-b}{}^2+D_F\para{x_{k+1},x_k}+\gamma_k^2\frac{L_k}{2}\norm{\nsk-x_k}{}^2 +\gamma_k\para{ \ip{\Enk, \xs -\nsk}{} +\snk}.
	\end{eqnarray*}
	The previous inequality can be re-written, by trivial manipulations, as\fekn
	\begin{equation}\label{est}
	\begin{split}
	& (1-c\theta_{k+1})\Lkp{x_{k+1},\mu_{k+1}} -(1-c\theta_k)\Lk{x_k,\mu_k} + \frac{c}{2}\para{\norm{\mu_{k+1}-\mus}{}^2-\norm{\mu_k-\mus}{}^2}\\
	\leq & \ \ (1-c\theta_{k+1})\Lkp{x_{k+1},\mu_{k+1}}-(1-c\theta_{k})\Lk{x_{k+1},\mu_{k+1}}+\para{\theta_k -\frac{c\theta_k^2}{2}}\norm{Ax_{k+1}-b}{}^2\\
	& \ +\gamma_k \para{\Lk{\xs,\mu_k}-\Lk{x_k,\mu_k}}-c\theta_k \para{\Lk{x_{k+1},\mus}-\Lk{x_k,\mu_k}}-\frac{\rho_k\gamma_k}{2}\norm{Ax_{k}-b}{}^2 \\
	& \ +D_F\para{x_{k+1},x_k}+\gamma_k^2\frac{L_k}{2}\norm{\nsk-x_k}{}^2 +\gamma_k\para{ \ip{\Enk, \xs -\nsk}{} +\snk}\\
	=&\ \ c\para{\theta_k-\theta_{k+1}}\para{f+h+\scal{\mu_{k+1}}{A\cdot-b}}(x_{k+1})+\para{\para{1-c\theta_{k+1}}g^{\beta_{k+1}}-\para{1-c\theta_k}g^{\beta_k}}\para{Tx_{k+1}}\\
	& \ + \frac{1}{2}\para{\para{1-c\theta_{k+1}}\rho_{k+1}-\para{1-c\theta_{k}}\rho_{k}+2\theta_k -c\theta_k^2}\norm{Ax_{k+1}-b}{}^2+\gamma_k \para{\Lk{\xs,\mu_k}-\Lk{x_k,\mu_k}}\\
	& \ -c\theta_k \para{\Lk{x_{k+1},\mus}-\Lk{x_k,\mu_k}}-\frac{\rho_k\gamma_k}{2}\norm{Ax_{k}-b}{}^2 +D_F\para{x_{k+1},x_k}+\gamma_k^2\frac{L_k}{2}\norm{\nsk-x_k}{}^2\\
	& \ +\gamma_k\para{ \ip{\Enk, \xs -\nsk}{} +\snk}.
	\end{split}
\end{equation}
By \ref{ass:gkgkp1}, \ref{ass:thetagamma} and the assumption that $\gamconinf \geq 1$, we have $\theta_{k+1} \leq \gamconinf^{-1}\theta_k \leq \theta_k$. In view of \ref{ass:beta}, we also have $\beta_{k+1}\leq\beta_{k}$. In particular, $g^{\beta_{k}}\leq g^{\beta_{k+1}}\leq g$ pointwise. By Proposition \ref{moreau}\ref{moreauclaim3} and assumption \ref{ass:interior}, we are able to\fekn estimate the quantity
	\begin{equation*}
	\begin{split}	
	&\para{\para{1-c\theta_{k+1}}g^{\beta_{k+1}}-\para{1-c\theta_k}g^{\beta_k}}\para{Tx_{k+1}}\\
	&=\para{g^{\beta_{k+1}}-g^{\beta_k}}\para{Tx_{k+1}}+c\para{\theta_kg^{\beta_k}-\theta_{k+1}g^{\beta_{k+1}}}\para{Tx_{k+1}}\\
	&\leq \frac{1}{2}\para{\beta_k-\beta_{k+1}}\norm{\para{\partial g(Tx_{k+1})}^0}{}^2+c\para{\theta_kg^{\beta_k}-\theta_{k+1}g^{\beta_{k}}}\para{Tx_{k+1}}\\
	&\leq \frac{1}{2}\para{\beta_k-\beta_{k+1}}\norm{\para{\partial g(Tx_{k+1})}^0}{}^2+c\para{\theta_k-\theta_{k+1}}g(Tx_{k+1}).
	\end{split}
	\end{equation*}
	Then\fekn
	\begin{equation}\label{est2}
	\begin{split}
	& \ \ c\para{\theta_k-\theta_{k+1}}\para{f+h+\scal{\mu_{k+1}}{A\cdot-b}}(x_{k+1})+\para{\para{1-c\theta_{k+1}}g^{\beta_{k+1}}-\para{1-c\theta_k}g^{\beta_k}}\para{Tx_{k+1}}\\
	\leq &\ \ c\para{\theta_k-\theta_{k+1}}\LL{x_{k+1},\mu_{k+1}}+\frac{1}{2}\para{\beta_k-\beta_{k+1}}\norm{\para{\partial g(Tx_{k+1})}^0}{}^2.
	\end{split}
	\end{equation}
	
	Recall the definition of $r_k$ in \eqref{r}. Coming back to \eqref{est} and using \eqref{est2}, we obtain\fekn
	\begin{equation}\label{est3}
	\begin{split}
	r_{k+1}-r_k
	\ & \leq \ \frac{1}{2}\para{\para{1-\gamma_{k+1}}\rho_{k+1}-\para{1-\gamma_{k}}\rho_{k}+\frac{2}{c}\gamma_k -\frac{\gamma_k^2}{c}}\norm{Ax_{k+1}-b}{}^2+\gamma_k \para{\Lk{\xs,\mu_k}-\Lk{x_{k+1},\mus}}\\
	& \ \  -\frac{\rho_k\gamma_k}{2}\norm{Ax_{k}-b}{}^2 +\frac{\beta_k-\beta_{k+1}}{2}\norm{\para{\partial g(Tx_{k+1})}^0}{}^2+ \para{\gamma_k-\gamma_{k+1}}\LL{x_{k+1},\mu_{k+1}}\\
	& \ \ +D_F\para{x_{k+1},x_{k}}+\gamma_k^2\frac{L_k}{2}\norm{\nsk-x_k}{}^2 +\gamma_k\para{ \ip{\Enk, \xs -\nsk}{} +\snk}.
	\end{split}
	\end{equation}
Recall that, by feasibility of $\xs$ for the affine constraint, $\LL{\xs,\mu_k}=\LL{\xs,\mus}$ and thus\fekn
	\begin{equation*}
	\begin{split}
	\Lk{\xs,\mu_k}-\Lk{x_{k+1},\mus} & = \LL{\xs,\mus}-\LL{x_{k+1},\mus}+\para{g^{\beta_k}-g}(T\xs)+ \para{g-g^{\beta_{k}}}(Tx_{k+1})\\
	&\quad\quad-\frac{\rho_k}{2}\norm{Ax_{k+1}-b}{}^2\\
	&=\LL{\xs,\mus}-\LL{x_{k},\mus} + \LL{x_{k},\mus}-\LL{x_{k+1},\mus} \\
	&\quad\quad\para{g^{\beta_k}-g}(T\xs)+ \para{g-g^{\beta_{k}}}(Tx_{k+1})-\frac{\rho_k}{2}\norm{Ax_{k+1}-b}{}^2\\
	&\leq\LL{\xs,\mus}-\LL{x_{k},\mus} + \LL{x_{k},\mus}-\LL{x_{k+1},\mus} + \frac{\beta_{k}}{2}\norm{\para{\partial g(Tx_{k+1})}^0}{}^2\\
	&\quad\quad-\frac{\rho_k}{2}\norm{Ax_{k+1}-b}{}^2,
	\end{split}
	\end{equation*}
where in the inequality we have used the fact that $g^{\beta_k}\leq g$ pointwise and that, by Proposition \ref{moreau}\ref{moreauclaim4}\fekn
	$$\para{g-g^{\beta_{k}}}(Tx_{k+1}) \leq \frac{\beta_{k}}{2}\norm{\para{\partial g(Tx_{k+1})}^0}{}^2.$$ 

Substituting the above into \eqref{est3} we have\fekn
\nnewq{\label{est4}
	\begin{split}
	r_{k+1}-r_k
	\ & \leq \ \frac{1}{2}\para{\para{1-\gamma_{k+1}}\rho_{k+1}-\rho_{k}+\frac{2}{c}\gamma_k -\frac{\gamma_k^2}{c}}\norm{Ax_{k+1}-b}{}^2\\
	& \ \ +\gamma_k \para{\LL{\xs,\mus}-\LL{x_{k},\mus}} + \gamma_{k}\para{\LL{x_{k},\mus}-\LL{x_{k+1},\mus}}\\
	& \ \  -\frac{\rho_k\gamma_k}{2}\norm{Ax_{k}-b}{}^2 +\frac{\beta_k-\beta_{k+1}}{2}\norm{\para{\partial g(Tx_{k+1})}^0}{}^2 + \para{\gamma_k-\gamma_{k+1}}\LL{x_{k+1},\mu_{k+1}}\\
	& \ \  + \gamma_k\frac{\beta_{k}}{2}\norm{\para{\partial g(Tx_{k+1})}^0}{}^2 + D_F\para{x_{k+1},x_{k}}+\gamma_k^2\frac{L_k}{2}\norm{\nsk-x_k}{}^2\\
	& \ \ \  +\gamma_k\para{ \ip{\Enk, \xs -\nsk}{} +\snk}.
	\end{split}
}
\ \\
Now we take the expectation with respect to $\filtop_{k}=\filts_{k-1} = \sigma\para{x_0,\mu_0,\nsz,\ldots,\nskm}$, which will completely determine $x_k$ and $\mu_k$, and we are perform the following estimations.
\begin{itemize}
	\item  From \ref{ass:cond}, we have\fekn
	$$\para{1-\gamma_{k+1}}\rho_{k+1}-\rho_{k}+\frac{2}{c}\gamma_k -\frac{\gamma_k^2}{c}\leq\gamma_{k+1}.$$
	\item By assumption \ref{ass:interior}\fekn
	$$ \EX{\norm{\para{\partial g(Tx_{k+1})}^0}{}^2}{\filtop_{k}}\leq M^2.$$ 
	\item By \lemref{unifbound}\fekn $$\EX{\LL{x_{k+1},\mu_{k+1}}}{\filtop_{k}}\leq \tilde{M}.$$
	\item Recall that, by \ref{ass:f}, $f$ is $\para{F,\zeta}$-smooth and invoke \remref{fremark}, to get\fekn
	$$\EX{D_F\para{x_{k+1},x_{k}}}{\filtop_{k}}\leq\KFC\zeta\para{\gamma_k}.$$
	\item Since\fekn $\nsk$ and $x_{k}$ are both in $\C$, we have $$\EX{\norm{\nsk-x_k}{}}{\filtop_{k}}\leq \diam.$$
\end{itemize}

We have\fekn
\newq{
	\begin{split}
\EX{r_{k+1}}{\filtop_{k}} - r_k &\leq \frac{\gamma_{k+1}}{2}\EX{\norm{Ax_{k+1}-b}{}^2}{\filtop_{k}}+\gamma_k \para{\LL{\xs,\mus}-\LL{x_{k},\mus}}\\
& \ \ +\gamma_k\para{\LL{x_{k},\mus}-\EX{\LL{x_{k+1},\mus}}{\filtop_{k}}}-\frac{\rho_k\gamma_k}{2}\norm{Ax_{k}-b}{}^2 +\frac{\beta_k-\beta_{k+1}}{2}M^2\\
	& \ \  + \para{\gamma_k-\gamma_{k+1}}\tilde{M}+ \gamma_k\frac{\beta_{k}}{2}M^2+ \KFC\zeta\para{\gamma_k}+\gamma_k^2\frac{L_k}{2}\diam^2+\gamma_k\EX{ \ip{\Enk, \xs -\nsk}{} +\snk}{\filtop_{k}}.
	\end{split}
}
We can bound the inner product involving the error terms using the Cauchy-Schwartz inequality and the boundedness of $\C$. Applying \lemref{convopt} and regrouping terms with $\gamma_k^2$ we get\fekn
\newq{
	\begin{split}
\EX{r_{k+1}}{\filtop_{k}} - r_k &\leq \frac{\gamma_{k+1}}{2}\EX{\norm{Ax_{k+1}-b}{}^2}{\filtop_{k}}+\gamma_k \para{\LL{\xs,\mus}-\LL{x_{k},\mus}}-\frac{\rho_k\gamma_k}{2}\norm{Ax_{k}-b}{}^2\\
	& \ \  +\para{\beta_k-\beta_{k+1}}\frac{\gcon^2}{2}+ \para{\gamma_k-\gamma_{k+1}}\Mcstopt+ \gamma_k\beta_{k}\frac{\gcon^2}{2}+ \KFC\zeta\para{\gamma_k}+\gamma_k^2C_k\\
	& \ \ \ +\gamma_k\EX{ \diam\para{\norm{\Enk}{}} +\snk}{\filtop_{k}}.
	\end{split}
}
We conclude by trivial manipulations.
\end{proof}

We now proceed to prove the main theorem regarding optimality.
\begin{theorem}[Optimality]\label{opttheorem}
Suppose that \ref{ass:A1}-\ref{ass:coercive} and \ref{ass:P1}-\ref{ass:error} hold, with $\gamconinf \geq 1$. Let $\seq{x_k}$ be the sequence of primal iterates generated by Algorithm~\ref{alg:CGAL} and $(\xs,\mus)$ a saddle-point pair for the Lagrangian. Then, in addition to the results of Theorem~\ref{feastheorem}, the following holds
\begin{enumerate}[label=(\roman*)]
\item \label{convlag}Convergence of the Lagrangian:
	\begin{equation}\label{limLL}
	\begin{aligned}
		\lim_{k\to\infty} \LL{x_{k},\mus} = \LL{\xs,\mus} \Pas.
	\end{aligned}
	\end{equation}
\item \label{wkcluster}Every weak cluster point $\bar{x}$ of $\seq{x_k}$ is a solution of the primal problem \eqref{PProb}, and $\seq{\mu_k}$ converges weakly to $\bar{\mu}$ a solution of the dual problem \eqref{DProb}, i.e., $(\bar{x},\bar{\mu})$ is a saddle point of $\mc{L}$ $\Pas$.
\item \label{pwrate}Pointwise rate: 
	\begin{equation}\label{pointrateopt}
	\begin{split}
		\forall k\in\N, \inf_{0 \leq i \leq k} \LL{x_{i},\mus}-\LL{\xs,\mus} = O\para{\frac{1}{\Gamma_{k}}} \Pas \qandq \\
		\text{$\exists$ a subsequence $\subseq{x_{k_j}}$ s.t. }\forall j\in\N, \LL{x_{k_j+1},\mus}-\LL{\xs,\mus} \leq \frac{1}{\Gamma_{k_j}} \Pas.
	\end{split}
	\end{equation}
\item \label{ergrate}Ergodic rate: for each $k\in\N$, let $\avx_k \eqdef \sum_{i=0}^{k} \gamma_ix_{i+1}/\Gamma_k$. Then\fekn
\begin{equation}\label{ergrateopt}
	\begin{split}
	\LL{\avx_{k},\mus}-\LL{\xs,\mus} = O\para{\frac{1}{\Gamma_{k}}} \Pas.
	\end{split}
	\end{equation}
\item \label{unisol} If the problem \eqref{PProb} admits a unique solution $\xs$, then the primal-dual pair sequence $\seq{x_k,\mu_k}$ converges weakly$\Pas$ to a saddle point $\para{\xs,\mus}$. Moreover, if $\Phi$ is uniformly convex on $\C$ with modulus of convexity $\psi:\R_+\to[0,\infty]$, then $\seq{x_k}$ converges strongly$\Pas$ to $\xs$ at the ergodic rate\fekn
	\newq{
	\psi\para{\norm{\avx_k-\xs}{}}=O\para{\frac{1}{\Gamma_k}}\Pas.
	}
\end{enumerate}
\end{theorem}
\begin{proof}
As in the proof of \thmref{feastheorem}, our goal is to first apply \lemref{combstoch} and then apply \lemref{barty}. By \lemref{optlem} we have, using the same notation\fekn
\newq{
	\begin{split}
\EX{r_{k+1}}{\filtop_{k}} - r_k  &\leq -\gamma_k \para{\LL{x_{k},\mus}-\LL{\xs,\mus} + \frac{\rho_k}{2}\norm{Ax_k-b}{}^2} + \frac{\gamma_{k+1}}{2}\EX{\norm{Ax_{k+1}-b}{}^2}{\filtop_{k}}\\
	& \ \  +\para{\beta_k-\beta_{k+1}}\frac{\gcon^2}{2}+ \para{\gamma_k-\gamma_{k+1}}\Mcstopt+ \gamma_k\beta_{k}\frac{\gcon^2}{2}+ \KFC\zeta\para{\gamma_k}+\gamma_k^2C_k\\
	& \ \ \ +\diam\gamma_k\EX{\norm{\Enk}{}}{\filtop_{k}} +\gamma_k\EX{\snk}{\filtop_{k}}.
	\end{split}
}
Let\fekn $a_k = \gamma_k \para{\LL{x_{k},\mus}-\LL{\xs,\mus} + \frac{\rho_k}{2}\norm{Ax_k-b}{}^2}$ and denote what remains on the r.h.s. by $z_k$. Then, to apply \lemref{combstoch}, we must show $\seq{z_k}\in\ell^1_+\para{\Filt}$. The first term, $\gamma_{k+1}\EX{\norm{Ax_{k+1}-b}{}^2}{\Filtop_k}$, is in $\ell^1_+\para{\Filt}$ by \ref{feastheorem}. The terms $\para{\beta_k-\beta_{k+1}}\frac{\gcon^2}{2}$ and $\para{\gamma_k-\gamma_{k+1}}\Mcstopt$ are bounded and telescopic, hence in $\ell^1_+$. The terms $\gamma_k\beta_k \frac{\gcon^2}{2}$ and $\KFC\zeta\para{\gamma_k}$ are in $\ell^1_+$ by \ref{ass:P1}. Recalling the definition of $C_k$, we have\fekn
\newq{
\gamma_k^2C_k 	&= \gamma_k^2\para{\frac{L_k}{2}\diam^2 + \diam\para{\gcon\norrm{T}+D+L_h+\norrm{\mus}\ \norrm{A}}}\\
				&= \para{\frac{\diam^2\norm{T}{}^2}{2}}\frac{\gamma_k^2}{\beta_k} + \para{\frac{\diam^2\norm{A}{}^2\rho_k}{2}+ \diam\para{\gcon\norrm{T}+D+L_h+\norrm{\mus}\ \norrm{A}}}\gamma_k^2\\
				&\leq \para{\frac{\diam^2\norm{T}{}^2}{2}}\frac{\gamma_k^2}{\beta_k} + \para{\frac{\diam^2\norm{A}{}^2\rhosup}{2}+ \diam\para{\gcon\norrm{T}+D+L_h+\norrm{\mus}\ \norrm{A}}}\gamma_k^2
}
which is in $\ell^1_+$ by \ref{ass:P1} and \ref{ass:beta}. The remaining terms,
\newq{
\diam\gamma_k\EX{\norm{\Enk}{}}{\filtop_{k}} +\gamma_k\EX{\snk}{\filtop_{k}},
}
 coming from the inexactness of the algorithm, are in $\ell^1_+\para{\Filt}$ by \ref{ass:error}. Thus, the r.h.s. belongs to $\ell^1_+\para{\Filt}$ and so by \lemref{combstoch} we have,
\newq{
a_k = \gamma_k \para{\LL{x_{k},\mus}-\LL{\xs,\mus} + \frac{\rho_k}{2}\norm{Ax_k-b}{}^2}\in\ell^1_+\para{\Filt}\Pas.
}
The first claim \ref{convlag} follows by applying \lemref{barty}, the conditions of which are satisfied directly from \lemref{convfeas} and \lemref{convopt}. The following three claims, \ref{wkcluster}, \ref{pwrate}, and \ref{ergrate}, all follow from \cite[Theorem~4.2]{silveti}. The final claim, \ref{unisol}, follows from \cite[Corollary~19]{silveti}.
\end{proof}

\ifdefined\COMPLETE
\else
\end{document}
\fi

\ifdefined\COMPLETE
\else
\documentclass[12pt]{article}
\input{package_header} 
\begin{document}
\fi

\section{Stochastic Examples}\label{sec:exam}
We examine the problem of risk minimization using two different ways to inexactly calculate the gradient with stochastic noise to demonstrate that the assumptions on the error can be satisfied in order to apply \icgalp. 

Consider the following,
\begin{equation}
\label{BatProb}\tag{$\mathrsfs{P}_1$}
\min\limits_{\sst{x\in \C \subset \HH\\ Ax = b}} f\para{x}\sbrac{\eqdef \TEX{L\para{x,\eta}}}
\end{equation}
where $L\para{\cdot,\eta}$ is differentiable for every $\eta$, and $\eta$ is a random variable. 

We will impose the following assumptions, or a subset of them depending on the context:
\be[label=(E.\arabic*)]
\item It holds $\nabla_x f\para{x} = \TEX{\nabla_x L\para{x,\eta}}\Pae$\label{exass1}
\item For all $\eta$, the function $L\para{\cdot, \eta}$ is $\omega$-smooth (see \defref{def:omega}) with $\omega$ nondecreasing\label{exass2}
\item The function $f$ is $\omega$-smooth with $\omega$ nondecreasing\label{exass3}
\item The function $f$ is H\"{o}lder smooth with constant $C_f$ and exponent $\tau$.\label{exass4}
\ee

Notice that \ref{exass4}$\implies$\ref{exass3}. For the sake of clarity, we analyze only the case where\fekn $\Enk \equiv \fnk$ with $\fnk = \nfg_k-\nabla f\para{x_k}$ and $\nfg_k$ is our inexact computation of $\nabla f\para{x_k}$, to be defined in the following subsections. 

\begin{remark}\label{Esamp}
With the above choice for $\Enk$, the terms in $\grEk{x_k,\mu_k}$ coming from the augmented Lagrangian are computed exactly, however our analysis extends to the case where $\nabla_x \para{\frac{\rho_k}{2}\norm{Ax_k-b}{}^2} = \rho_k A^*\para{Ax_k-b}$ is computed inexactly as well, as this function is always Lipschitz-continuous. We demonstrate this alternative choice in \secref{sec:numexp} by sampling the components $\rho_k \comp{A^*\para{Ax_k-b}}{i}$ in the numerical experiments.
\end{remark}

\subsection{Risk minimization with increasing batch size}
Consider \eqref{BatProb} and define\fekn  
\newq{
\nfg_k\eqdef \frac{1}{n\para{k}}\sum\limits_{i=1}^{n\para{k}}\nabla_x L\para{x_k,\eta_i}
}
 where $n\para{k}$ is the number of samples to be taken at iteration $k$. We assume that each $\eta_i$ is i.i.d., according to some fixed distribution, and that $n$ is a function of $k$, i.e., the number of samples taken to estimate the expectation is dependent on the iteration number itself. 

\begin{lemma}\label{lem:ngrow}
Under assumptions \ref{exass1} and \ref{exass2}, denote
\newq{
C=2\para{\omega\para{d_\C}^2 + \EX{\norm{\nabla L\para{\xs,\eta}}{}^2}{\filts_k}}
}
where $\xs$ is a solution to \eqref{BatProb} and\fekn $\filts_k=\sigma\para{x_0,\mu_0,\what{s}_0,\ldots,\nsk}$ as before.
Then, for each $k\in\N$, the following holds,
\newq{
\EX{\norm{\fnkp}{}}{\filts_k}\leq \sqrt{\frac{C}{n\para{k+1}}}.
}
\end{lemma}
\begin{proof}
By Jensen's inequality\fekn
\newq{
\EX{\norm{\fnkp}{}}{\filts_k}^2 	&\leq \EX{\norm{\fnkp}{}^2}{\filts_k} = \EX{\norm{\nabla f\para{x_{k+1}} - \nfg_{k+1}}{}^2}{\filts_k}.
}
Then, since $\nfg_{k+1}$ is an unbiased estimator for $\nabla f\para{x_{k+1}}$, we have\fekn
\newq{
\EX{\norm{\nabla f\para{x_{k+1}} - \nfg_{k+1}}{}^2}{\filts_k}								&= \EX{\norm{\TEX{\nfg_{k+1}} - \nfg_{k+1}}{}^2}{\filts_k}\\
								&=\VAR{\nfg_{k+1}}{\filts_k}\\
								&=\VAR{\frac{1}{n\para{k+1}} \sum\limits_{i=1}^{n\para{k+1}} \nabla L\para{x_{k+1},\eta_i}}{\filts_k}\\
								&=\frac{1}{n\para{k+1}}\VAR{\nabla L\para{x_{k+1},\eta}}{\filts_k},
}
where the last equality follows from the independence and identical distribution of $\eta_i$. Applying the definition of conditional variance yields\fekn
\newq{
\frac{1}{n\para{k+1}}\VAR{\nabla L\para{x_{k+1},\eta}}{\filts_k}								&=\frac{1}{n\para{k+1}}\para{\EX{\norm{\nabla L\para{x_{k+1},\eta}}{}^2}{\filts_k} - \norm{\EX{\nabla L\para{x_{k+1},\eta}}{\filts_k}}{}^2}\\
								&\leq \frac{1}{n\para{k+1}}\EX{\norm{\nabla L\para{x_{k+1},\eta}}{}^2}{\filts_k}.
}
We again use Jensen's inequality, then $\omega$-smoothness, and finally the fact that $\omega$ is nondecreasing together with the fact that $x_{k+1}$ and $\xs$ are both in $\C$ to find\fekn
\newq{
\frac{1}{n\para{k+1}}\EX{\norm{\nabla L\para{x_{k+1},\eta}}{}^2}{\filts_k}								&\leq \frac{2}{n\para{k+1}}\para{\EX{\norm{\nabla L\para{x_{k+1},\eta}-\nabla L\para{\xs,\eta}}{}^2}{\filts_k} \right.\\
											&\left.\quad\quad+ \EX{\norm{\nabla L\para{\xs,\eta}}{}^2}{\filts_k}}\\
								&\leq \frac{2}{n\para{k+1}}\para{\EX{\omega\para{\norm{x_{k+1}-\xs}{}}^2}{\filts_k} + \EX{\norm{\nabla L\para{\xs,\eta}}{}^2}{\filts_k}}\\
								&\leq \frac{2}{n\para{k+1}}\para{\omega\para{d_\C}^2 + \EX{\norm{\nabla L\para{\xs,\eta}}{}^2}{\filts_k}}\\
								&= \frac{C}{n\para{k+1}}.
}
The above shows that\fekn $\EX{\norm{\fnkp}{}}{\filts_k}^2\leq \frac{C}{n\para{k+1}}$ and so $\EX{\norm{\fnkp}{}}{\filts_k}\leq \sqrt{\frac{C}{n\para{k+1}}}$ as desired.
\end{proof}
\begin{proposition}\label{prop:ngrow}
Under \ref{exass1} and \ref{exass2}, assume that the number of samples $n\para{k}$ at iteration $k$ is lower bounded by $\para{\frac{\gamma_k}{\zeta\para{\gamma_k}}}^{2}$, i.e. for some $\alpha >0$, $n\para{k} \geq \alpha \para{\frac{\gamma_k}{\zeta\para{\gamma_k}}}^{2}$. Then, the summability of the error in \ref{ass:error} is satisfied; namely,
\newq{
\gamma_{k+1}\EX{\norm{\fnkp}{}}{\filts_k} \in\ell^1\para{\Filts}.
}

\end{proposition}
\begin{proof}
By \lemref{lem:ngrow} we have\fekn
\newq{
\gamma_{k+1}\EX{\norm{\fnkp}{}}{\filts_k} \leq \gamma_{k+1} \sqrt{\frac{C}{n\para{k+1}}} \leq \sqrt{\frac{C}{\alpha }}\zeta\para{\gamma_{k+1}}.
}
The summability of $\zeta\para{\gamma_{k+1}}$ is given by \ref{ass:P1} and thus $\gamma_{k=1}\EX{\norm{\fnkp}{}}{\filts_k}\in\ell^1\para{\Filts}$
\end{proof}
\begin{remark}
The lower bound $n\para{k} \geq \alpha \para{\frac{\gamma_k}{\zeta\para{\gamma_k}}}^2$ is sufficient but not necessary; one can alternatively choose $n\para{k}$ to be lower bounded by $\alpha\para{\frac{\beta_k}{\gamma_k}}^2$ or $\alpha\para{\frac{1}{\beta_k}}^2$ and, due to \ref{ass:P1}, the result will still hold.
\end{remark}

\subsection{Risk minimization with variance reduction}

We reconsider \eqref{BatProb} as before but now with a different $\nfg$. We define a stochastic-averaged gradient, which will serve as a form of variance reduction, such that the number of samples at each iteration need not increase as in the previous subsection. For each $k\in\N$, let $\nu_k\in[0,1]$ and define
\nnewq{\label{nfgvr}
\nfg_k\eqdef \para{1-\nu_k}\nfg_{k-1} + \nu_k \nabla_x L\para{x_k,\eta_k}
}
 with $\nfg_{-1}=0$ and with each $\eta_i$ i.i.d.. We call $\nfg_k$ the stochastic average of sampled gradients with weight $\nu_k$. In this way, we are able to take a single gradient sample (or a larger fixed batch size) at each iteration, in contrast to the previous subsection.
\begin{lemma}\label{lem:stochavg}
Under \ref{exass1} and \ref{exass3}, denote, for each $k\in\N$,
\nnewq{\label{sigdef}
\sigma_k^2 \eqdef \EX{\norm{\nabla_x L \para{x_{k},\eta_k} - \nabla f\para{x_k}}{}^2}{\filts_k}
}
 and assume that $\exists \sigma > 0$ such that $\sup_k \sigma_k^2 = \sigma^2 <\infty$. Then, for each $k\in\N$, the following inequality holds,
\newq{
\EX{\norm{\fnkp}{}^2}{\filts_{k}} \leq \para{1-\frac{\nu_{k+1}}{2}}\norm{\fnk}{}^2 + \nu_{k+1}^2\sigma^2 +2\frac{\omega\para{d_\C \gamma_k}^2}{\nu_{k+1}}.
}
\end{lemma}
\begin{proof}
The proof of this theorem is inspired by a similar construction found in \cite[Lemma~2]{SARAH}. By definition of $\fnkp$ and $\nfg_{k+1}$, we have, for all $k\in\N$,
\newq{
\norm{\fnkp}{}^2 = \norm{\nfg_{k+1}-\nabla f \para{x_{k+1}}}{}^2 = \norm{\para{1-\nu_{k+1}}\nfg_{k} + \nu_{k+1}\nabla_x L\para{x_{k+1},\eta_{k+1}} - \nabla f\para{x_{k+1}}}{}^2.
}
We add and subtract $\para{1-\nu_{k+1}}\nabla f\para{x_k}$ to get,
\newq{
\norm{\fnkp}{}^2 = \norm{\para{1-\nu_{k+1}}\fnk + \nu_{k+1}\para{\nabla_x L\para{x_{k+1},\eta_{k+1}} - \nabla f\para{x_{k+1}}} + \para{1-\nu_{k+1}}\para{\nabla f\para{x_k}-\nabla f\para{x_{k+1}}}}{}^2.
}
Applying the pythagoreas identity then gives,
\newq{
\norm{\fnkp}{}^2	&= \para{1-\nu_{k+1}}^2 \norm{\fnk}{}^2 + \nu_{k+1}^2\norm{\nabla_x L\para{x_{k+1},\eta_{k+1}} - \nabla f\para{x_{k+1}}}{}^2\\
					&\quad\quad + \para{1-\nu_{k+1}}^2\norm{\nabla f\para{x_k}-\nabla f\para{x_{k+1}}}{}^2\\
					& \quad\quad + 2\ip{\para{1-\nu_{k+1}}\para{\fnk + \nabla f\para{x_k}-\nabla f\para{x_{k+1}}}, \nu_{k+1}\para{\nabla_x L\para{x_{k+1},\eta_{k+1}} - \nabla f\para{x_{k+1}}}}{}\\
					& \quad\quad + 2\ip{\para{1-\nu_{k+1}}\fnk, \para{1-\nu_{k+1}}\para{\nabla f\para{x_k}-\nabla f\para{x_{k+1}}}}{}.
}
Using Young's inequality on the last inner product, we find,
\newq{
\norm{\fnkp}{}^2 &\leq \para{1-\nu_{k+1}}^2\norm{\fnk}{}^2 + \nu_{k+1}^2\norm{\nabla_x L\para{x_{k+1},\eta_{k+1}} - \nabla f\para{x_{k+1}}}{}^2\\
					& \quad\quad + \para{1-\nu_{k+1}}^2\norm{\nabla f\para{x_k}-\nabla f\para{x_{k+1}}}{}^2\\
					& \quad\quad + 2\ip{\para{1-\nu_{k+1}}\para{\fnk + \nabla f\para{x_k}-\nabla f\para{x_{k+1}}}, \nu_{k+1}\para{\nabla_x L\para{x_{k+1},\eta_{k+1}} - \nabla f\para{x_{k+1}}}}{}\\
					& \quad\quad + \frac{\nu_{k+1}}{2}\norm{\fnk}{}^2 + \frac{2}{\nu_{k+1}}\norm{\para{1-\nu_{k+1}}^2\para{\nabla f\para{x_k}-\nabla f\para{x_{k+1}}}}{}^2.
}
Notice that $1-\nu_{k+1}\leq 1$ and thus $\para{1-\nu_{k+1}}^2\leq 1-\nu_{k+1}$ for all $k\in\N$. This leads to
\newq{
\norm{\fnkp}{}^2 &\leq \para{1-\frac{\nu_{k+1}}{2}}\norm{\fnk}{}^2 + \nu_{k+1}^2\norm{\nabla_x L\para{x_{k+1},\eta_{k+1}} - \nabla f\para{x_{k+1}}}{}^2 + \norm{\nabla f\para{x_k}-\nabla f\para{x_{k+1}}}{}^2\\
					& \quad\quad + 2\ip{\para{1-\nu_{k+1}}\para{\fnk + \nabla f\para{x_k}-\nabla f\para{x_{k+1}}}, \nu_{k+1}\para{\nabla_x L\para{x_{k+1},\eta_{k+1}} - \nabla f\para{x_{k+1}}}}{}\\
					& \quad\quad + \frac{2\para{1-\nu_{k+1}}}{\nu_{k+1}}\norm{\para{\nabla f\para{x_k}-\nabla f\para{x_{k+1}}}}{}^2\\
					&\leq \para{1-\frac{\nu_{k+1}}{2}}\norm{\fnk}{}^2 + \nu_{k+1}^2\norm{\nabla_x L\para{x_{k+1},\eta_{k+1}} - \nabla f\para{x_{k+1}}}{}^2 + \para{\frac{2}{\nu_{k+1}}}\norm{\nabla f\para{x_k}-\nabla f\para{x_{k+1}}}{}^2\\
					& \quad\quad + 2\ip{\para{1-\nu_{k+1}}\para{\fnk + \nabla f\para{x_k}-\nabla f\para{x_{k+1}}}, \nu_{k+1}\para{\nabla_x L\para{x_{k+1},\eta_{k+1}} - \nabla f\para{x_{k+1}}}}{}.
}
Recall that, by \ref{exass3}, $f$ is $\omega$-smooth with $\omega$ is nondecreasing. Furthermore, using the fact that $x_{k+1}=x_k-\gamma_{k}\para{x_k-\nsk}$, we find
\newq{
\norm{\fnkp}{}^2					&\leq \para{1-\frac{\nu_{k+1}}{2}}\norm{\fnk}{}^2 + \nu_{k+1}^2\norm{\nabla_x L\para{x_{k+1},\eta_{k+1}} - \nabla f\para{x_{k+1}}}{}^2 + \para{\frac{2}{\nu_{k+1}}}\omega\para{\norm{x_k-x_{k+1}}{}}^2\\
					& \quad\quad + 2\ip{\para{1-\nu_{k+1}}\para{\fnk + \nabla f\para{x_k}-\nabla f\para{x_{k+1}}}, \nu_{k+1}\para{\nabla_x L\para{x_{k+1},\eta_{k+1}} - \nabla f\para{x_{k+1}}}}{}\\
					&\leq \para{1-\frac{\nu_{k+1}}{2}}\norm{\fnk}{}^2 + \nu_{k+1}^2\norm{\nabla_x L\para{x_{k+1},\eta_{k+1}} - \nabla f\para{x_{k+1}}}{}^2 + \para{\frac{2}{\nu_{k+1}}}\omega\para{d_\C\gamma_{k}}^2\\
					& \quad\quad + 2\ip{\para{1-\nu_{k+1}}\para{\fnk + \nabla f\para{x_k}-\nabla f\para{x_{k+1}}}, \nu_{k+1}\para{\nabla_x L\para{x_{k+1},\eta_{k+1}} - \nabla f\para{x_{k+1}}}}{}
}
We take the expectation on both sides, recalling the definition of $\sigma_k$ (see \eqref{sigdef}), $\sigma$, and that
\newq{
\EX{\nabla_x L\para{x_k,\eta_k}}{\filts_{k-1}} = \nabla f\para{x_k},
}
to find,
\newq{
\EX{\norm{\fnkp}{}^2}{\filts_k} 					&\leq \para{1-\frac{\nu_{k+1}}{2}}\norm{\fnk}{}^2 + \nu_{k+1}^2\sigma^2 + \para{\frac{2}{\nu_{k+1}}}\omega\para{d_\C\gamma_{k}}^2.
}
\end{proof}

In the following proposition, we analyze a particular case of parameter choices under the assumption \ref{exass4} of H\"{o}lder smoothness of $f$, i.e. $\exists C_f,\tau >0$ such that $\omega:t\to C_ft^\tau$.

\begin{proposition}\label{prop:stochavg}
	Under \ref{exass1} and \ref{exass4}\fekn let $\nfg_k$ be defined as in \eqref{nfgvr} with weight $\nu_k=\gamma_k^\alpha$ for some $\alpha \in]0,\tau[$. If the following conditions on the sequence $\seq{\gamma_k}$ hold,
	\nnewq{\label{gammaconditions1}
		\seq{\gamma_k^{1+\min\brac{\frac{\alpha}{2}, \tau-\alpha}}}\in\ell^1,
	}
	and, for $k$ sufficiently large,
	\nnewq{\label{gammaconditions2}\frac{\gamma_k}{\gamma_{k+1}} \leq 1+o\para{\gamma_k^\alpha},
	}
	then the summability condition in \ref{ass:error} is satisfied; namely,
	\newq{
		\gamma_{k+1}\EX{\norm{\fnkp}{}}{\filts_k}\in\ell^1\para{\Filts}.
	}
\end{proposition}
\begin{proof}
	Since \ref{exass4}$\implies$\ref{exass3}, the assumptions \ref{exass1} and \ref{exass3} are satisfied and \lemref{lem:stochavg} gives, for all $k\in\N$,
	\newq{
		\EX{\norm{\fnkp}{}^2}{\filts_{k}} \leq \para{1-\frac{\gamma_{k+1}^\alpha}{2}}\norm{\fnk}{}^2 +\sigma^2\gamma_{k+1}^{2\alpha} + \frac{2C_f^2d_\C^{2\tau}\gamma_k^{2\tau}}{\gamma_{k+1}^\alpha}.
	}
	By \ref{ass:gkgkp1} we have, for all $k\in\N$, $\gamma_k\leq \gamcon\gamma_{k+1}$. It follows that\fekn
	\newq{
		\EX{\norm{\fnkp}{}^2}{\filts_{k}} \leq \para{1-\frac{\gamma_{k+1}^\alpha}{2}}\norm{\fnk}{}^2 +\sigma^2\gamma_{k+1}^{2\alpha} + 2\gamcon^{2\tau} C_f^2d_\C^{2\tau}\gamma_{k+1}^{2\tau-\alpha}.
	}
	Consolidating higher order terms gives\fekn
	\newq{
		\EX{\norm{\fnkp}{}^2}{\filts_{k}} \leq \para{1-\frac{\gamma_{k+1}^\alpha}{2}}\norm{\fnk}{}^2 +\para{\sigma^2 + 2\gamcon^{2\tau} C_f^2d_\C^{2\tau}}\gamma_{k+1}^{\min\brac{2\alpha, 2\tau-\alpha}}.
	}
	Since $\alpha <\tau\leq 1$ by \ref{gammaconditions1}, it holds that $\alpha < \min\brac{1, 2\tau-\alpha}$, and the first condition of \lemref{lem:gamest} is satisfied. Additionally, by \eqref{gammaconditions2}, we have that the second condition, \eqref{eq:gamass}, of \lemref{lem:gamest} is satisfied as well and we can apply \lemref{lem:gamest} with
	\newq{
		u_k = \norm{\fnk}{}^2,\quad c = \frac{1}{2},\quad s=\alpha,\quad d = \para{\sigma^2 + 2\gamcon^{2\tau} C_f^2d_\C^{2\tau}},\qandq t=\min\brac{2\alpha, 2\tau-\alpha},
	}
	to find, for $k$ sufficiently large,
	\newq{
		\EX{\norm{\fnkp}{}^2}{\filts_k} \leq 2\tilde C \gamma_{k+1}^{\min\brac{\alpha, 2\para{\tau-\alpha}}} + o\para{\gamma_{k+1}^{\min\brac{\alpha, 2\para{\tau-\alpha}}}}
	}
	and, by extension, for $k$ sufficiently large,
	\newq{
		\EX{\norm{\fnkp}{}}{\filts_k} \leq \sqrt{2\tilde C} \gamma_{k+1}^{\min\brac{\frac{\alpha}{2}, \tau-\alpha}} + o\para{\gamma_{k+1}^{\min\brac{\frac{\alpha}{2}, \tau-\alpha}}}.
	}
	Then, for $k$ sufficiently large,
	\newq{
		\gamma_{k+1}\EX{\norm{\fnkp}{}}{\filts_k} &\leq \gamma_{k+1}\para{\sqrt{2\tilde C} \gamma_{k+1}^{\min\brac{\frac{\alpha}{2}, \tau-\alpha}} + o\para{\gamma_{k+1}^{\min\brac{\frac{\alpha}{2}, \tau-\alpha}}}}\\
		&\leq \sqrt{2\tilde C}\gamma_{k+1}^{1+\min\brac{\frac{\alpha}{2}, \tau-\alpha}} + o\para{\gamma_{k+1}^{1+\min\brac{\frac{\alpha}{2}, \tau-\alpha}}}.
	}
	Under the assumptions \ref{gammaconditions1} we have $\gamma_k^{1+\min\brac{\frac{\alpha}{2}, \tau-\alpha}}\in\ell^1$ and thus the summability condition of \ref{ass:error} is satisfied.
\end{proof}
\begin{example}\label{examplegamma}
	The condition \eqref{gammaconditions1}
	in \propref{prop:stochavg} can be satisfied, for example, by taking $\gamma_k = \frac{1}{\para{k+1}^{1-b}}$. In this case, the condition \eqref{gammaconditions1} reduces to picking $b$ such that the following holds,
	\newq{
		\para{1-b}\para{1+\min\brac{\frac{\alpha}{2},\tau-\alpha}} > 1.
	}
	Rearranging, we find that this is equivalent to, 
	\nnewq{\label{bcondition}
		b < 1-\para{1+\min\brac{\frac{\alpha}{2}, \tau-\alpha}}^{-1}.
	}
	The condition \eqref{gammaconditions2} in \propref{prop:stochavg} can be satisfied under this choice of $\gamma_k$ as well. We have,
	\newq{
		\frac{\gamma_k}{\gamma_{k+1}} = \para{\frac{k+2}{k+1}}^{1-b} = \para{1+\frac{1}{k+1}}^{1-b}\approx 1 + \frac{1-b}{k+1} = 1 + o\para{\gamma_k^\epsilon}
	}
	for any $0<\epsilon<1$, for $k$ sufficiently large.
	
	Recall that the predicted convergence rates for the ergodic iterates $\avx_k$ given by \thmref{feastheorem} and \thmref{opttheorem} under this choice of step size are,
	\newq{
		\norrm{A\avx_k-b} = O\para{\frac{1}{\sqrt{\Gamma_{k}}}} \Pas\quad\mbox{and}\quad
		\LL{\avx_{k},\mus}-\LL{\xs,\mus} = O\para{\frac{1}{\Gamma_{k}}} \Pas,
	}
	where $\Gamma_k = \sum\limits_{i=0}^k\gamma_i=\sum\limits_{i=0}^k\frac{1}{\para{i+1}^{1-b}}$. Thus, choosing $b$ to be as large as possible is desired. For a given value of $\tau$ corresponding to the H\"{o}lder exponent of the gradient, the best choice for $\alpha$ is $\frac{2}{3}\tau$. If the problem is Lipschitz-smooth, then $\tau=1$ and we get $\alpha=\frac{2}{3}$.
	
	Notice that the choice of $\alpha$ does not directly affect the predicted rates of convergence, which now depend only on the constant $b$. However, the choice of $\alpha$ dictates the possible choices for $b$ which satisfy the assumptions and thus, indirectly, the rates of convergence as well. In the Lipschitz-smooth case, choosing $\alpha=\frac{2}{3}$ leads one to pick $b<1-\para{4/3}^{-1} = \frac{1}{4}$
\end{example}

\section{Sweeping}\label{sec:ex2}
We now consider an example in which the errors in the computation of $\nabla f$ are deterministic; a finite sum minimization problem,
\begin{equation}
\label{SweepProb}\tag{$\mathrsfs{P}_2$}
\min\limits_{\sst{x\in \C \subset \HH\\ Ax = b}} \frac{1}{n}\sum\limits_{i=1}^n f_i\para{x}
\end{equation}
where $n>1$ is fixed. We assume that:
\be[label=(F.\arabic*)]
\item $f_i$ is $\omega$-smooth (see \defref{def:omega}) for $1\leq i\leq n$ with $\omega$ nondecreasing\label{fexass1}
\item $\seq{\gamma_k}$ a nonincreasing sequence.\label{fexass2}
\ee

As in the previous section, \secref{sec:exam}, we examine only the case where\fekn $\Enk\equiv\fnk=\nabla f\para{x_k}-\nfg_k$, with $\nfg_k$ to be defined below, although our analysis is straightforward to adapt to the more general case where one computes $\rho_k A^*\para{Ax_k-b}$ inexactly as well, at the expense of brevity (see \remref{Esamp}). We will sweep, or cycle, through the functions $f_i$, taking the gradient of a single one at each iteration and recursively averaging with the past gradients. For notation, fixed $n$, we take $\modd{k} \eqdef \para{k \mod n}$ with the convention that $\modd{n} \eqdef n$. We define the inexact gradient in the following way,
\newq{
\nfg_k \eqdef \frac{1}{n}\sum\limits_{i=1}^k\nabla f_{i}\para{x_i} \quad\para{\forall k \leq n}
}
and
\newq{
	\nfg_k \eqdef \nfg_{k-1} + \frac{1}{n}\para{\nabla f_{\modd{k}}\para{x_k} - \nabla f_{\modd{k}}\para{x_{k-n}}}\quad\para{\forall k \geq n+1}.
}
For $k \geq n+1$ it can also be written in closed form as,
\newq{
\nfg_k = \frac{1}{n}\para{\sum\limits_{i=1}^{\modd{k}}\nabla f_i\para{x_{i+k -\modd{k}}} - \sum\limits_{i= \modd{k}+1}^n \nabla f_i\para{x_{ i+k-n-\modd{k}}} }.
}
\begin{lemma}\label{lem:sweep}
Let $C= \frac{1}{n}\para{n\para{n-1} + \para{n-1}\para{2n-1}}$. Under \ref{fexass1} and \ref{fexass2}, we then have,  for all $k\geq 2n-1$, the following,
\newq{
\norm{\fnkp}{} \leq C\omega\para{\gamma_{k+2-2n}d_\C}.
}
\end{lemma}
\begin{proof}
Using the definition of $\fnkp$ for $k\geq2n-1\geq n+1$, we have
\newq{
\norm{\fnkp}{}	&= \norm{\nabla f\para{x_{k+1}} - \nfg_{k+1}}{}\\
								&= \frac{1}{n}\norm{\para{\sum\limits_{i=1}^{\modd{k+1}} \nabla f_i\para{x_{k+1}} - \nabla f_i\para{x_{i+k+1 -\modd{k+1}}}}  \right.\\
								&\quad\quad \left. + \para{\sum\limits_{i= \modd{k+1}+1}^n \nabla f_i\para{x_{k+1}} - \nabla f_i\para{x_{i+k+1-n-\modd{k+1}}}}}{}.
}
Then, we apply the triangle inequality and $\omega$-smoothness of $f_i$ assumed in \ref{fexass1},
\newq{
\norm{\fnkp}{}								&\leq \frac{1}{n}\para{\sum\limits_{i=1}^{\modd{k+1}} \norm{\nabla f_i\para{x_{k+1}} - \nabla f_i\para{x_{i+k+1-\modd{k+1}}}}{} \right.\\
								&\quad\quad \left. + \sum\limits_{i=\modd{k+1}+1}^n \norm{\nabla f_i\para{x_{k+1}} - \nabla f_i\para{x_{i+k+1-n-\modd{k+1}}}}{}}\\
								&\leq \frac{1}{n}\para{\sum\limits_{i=1}^{\modd{k+1}} \omega\para{\norm{x_{k+1} - x_{i+k+1-\modd{k+1}}}{}} \right.\\
								&\quad\quad \left. + \sum\limits_{i=\modd{k+1}+1}^n \omega\para{\norm{x_{k+1} - x_{i+k+1-n-\modd{k+1}}}{}}}.
}
Now we add and subtract the iterates in between $x_{k+1}$ and $x_{i+k+1-\modd{k+1}}$ then use the definition $x_{k+1}=x_k+\gamma_k\para{\nsk-x_k}$ and the fact that, for all $k\in\N$, $\nsk$ and $x_k$ are in $\C$,
\newq{
\norm{\fnkp}{}								&\leq \frac{1}{n}\para{\sum\limits_{i=1}^{\modd{k+1}}\sum\limits_{j=1}^{ \ \ \ \modd{k+1}-i} \omega\para{\norm{x_{k+2-j} - x_{k+1-j}}{}} \right.\\
								&\quad\quad \left. + \sum\limits_{i=\modd{k+1}+1}^n \sum\limits_{j=1}^{\modd{k+1}-i+n}\omega\para{\norm{x_{k+2-j} - x_{k+1-j}}{}}}\\
								&\leq \frac{1}{n}\para{\sum\limits_{i=1}^{\modd{k+1}}\sum\limits_{j=1}^{ \ \ \ \modd{k+1}-i} \omega\para{\gamma_{k+1-j}d_\C} \right.\\
								&\quad\quad \left. + \sum\limits_{i=\modd{k+1}+1}^n \sum\limits_{j=1}^{\modd{k+1}-i+n}\omega\para{\gamma_{k+1-j}d_\C}}.
}
Recall that, by \ref{fexass2}, $\seq{\gamma_k}$ is nonincreasing, by \ref{fexass1}, $\omega$ is a nondecreasing function, and, for each $k\in\N$, $\modd{k}\leq n$. Then,
\newq{
\norm{\fnkp}{}								&\leq \frac{1}{n}\para{\sum\limits_{i=1}^{\modd{k+1}} \para{-i + \modd{k+1}}\omega\para{\gamma_{k+1+i - \modd{k+1}}d_\C}\right.\\
								&\quad\quad\left. + \sum\limits_{i=\modd{k+1}+1}^n \para{-i+n+\modd{k+1}}\omega\para{\gamma_{k+1+i-n-\modd{k+1}}d_\C}}\\
								&\leq \frac{1}{n}\para{ \modd{k+1}\para{-1 + \modd{k+1}}\omega\para{\gamma_{k+2 - \modd{k+1}}d_\C} \right.\\
								&\quad\quad \left. + \para{n-\modd{k+1}} \para{-1+n+\modd{k+1}}\omega\para{\gamma_{k+2-n-\modd{k+1}}d_\C}}\\
								&\leq \frac{1}{n}\para{ n\para{n-1}\omega\para{\gamma_{k+2 - n}d_\C} + \para{n-1}\para{2n-1}\omega\para{\gamma_{k+2-2n}d_\C}}\\
								&\leq \frac{1}{n} \para{n\para{n-1} + \para{n-1}\para{2n-1}}\omega\para{\gamma_{k+2-2n}d_\C}.
}
\end{proof}
\begin{proposition}\label{prop:sweep}
Under \ref{fexass1} and \ref{fexass2}, and assuming that $\seq{\gamma_{k}\omega\para{d_\C\gamma_{k}}}\in\ell^1$, the summability condition of \ref{ass:error} holds; namely,
\newq{
\gamma_{k+1}\norm{\fnkp}{}\in\ell^1.
}
\end{proposition}
\begin{proof}
By \lemref{lem:sweep}, we have, for all $k\geq2n-1$,
\newq{
\gamma_{k+1}\norm{\fnkp}{} \leq C\gamma_{k+1}\omega\para{d_\C\gamma_{k+2-2n}}\leq C \gamma_{k+2-2n}\omega\para{d_\C\gamma_{k+2-2n}}
}
where we have used the fact that $\seq{\gamma_k}$ is a nonincreasing sequence by \ref{fexass2}. Since $\seq{\gamma_k\omega\para{d_\C \gamma_k}}\in\ell^1$, the desired claim follows.
\end{proof}


\ifdefined\COMPLETE
\else
\end{document}
\fi

\section{Numerical Experiments}\label{sec:numexp}
We apply the sweeping method and the variance reduction method to solve the following projection problem,
\nnewq{\label{numprob}
\min\limits_{\sst{\norm{x}{1}\leq 1 \\ Ax=0}}\frac{1}{2n} \norm{x-y}{}^2,
}
where $x$ and $y$ are in $\R^n$. Notice that this problem fits both the risk minimization and the sweeping problem structures. By choosing $f_i\para{x} = \frac{1}{2}\para{x_i-y_i}^2$ we can rewrite the problem to apply the sweeping method of \secref{sec:ex2}. Alternatively, we can let $\eta$ be a random variable taking values in the set $\brac{1,\ldots,n}$ and write $L(x,\eta) = \frac{1}{2}\para{x_\eta - y_\eta}^2$ to cast the problem as risk minimization as in \secref{sec:exam}. In both of these cases, it is possible by our analysis to consider also sampling components of the components of the gradient term $\nabla_x \frac{\rho_k}{2}\norm{Ax_k}{}^2=\rho_k A^*Ax_k$.

The assumptions \ref{exass1} - \ref{exass4} and \ref{fexass1} all hold as the function $f$ is Lipschitz-smooth and the functions $L\para{\cdot,\eta}$ are all Lipschitz-smooth for every $\eta$ as well. The assumptions \eqref{ass:A1} to \eqref{ass:coerciveinfdim} all hold as $f$ is Lipschitz-smooth and has full domain.

For parameters, we take $\gamma_k=1/\para{k+1}^{1-b}$, $\rho_k\equiv \rho = 2^{2-b}+1$, $\theta_k=\gamma_k$. If we take $b<\frac{1}{2}$ then all the assumptions \ref{ass:P1} to \ref{ass:cond} are satisfied, as well as \ref{fexass2}. In particular, to satisfy \ref{ass:error} in the variance reduction case, we will take $b\in\brac{\frac{1}{4} - 0.15,\frac{1}{3} - 0.01}$. The weight $\nu_k$ in the variance reduction is chosen to be $\nu_k=\gamma_k^{\alpha}$ with $\alpha = 2/3$ since the problem is Lipschitz-smooth, i.e. the H\"{o}lder exponent is $\tau=1$. With this choice, the condition \eqref{gammaconditions1} in \propref{prop:stochavg} is satisfied as was discussed in Example \ref{examplegamma}.

\begin{figure}[hbp]
\centering
\includegraphics[width=\linewidth]{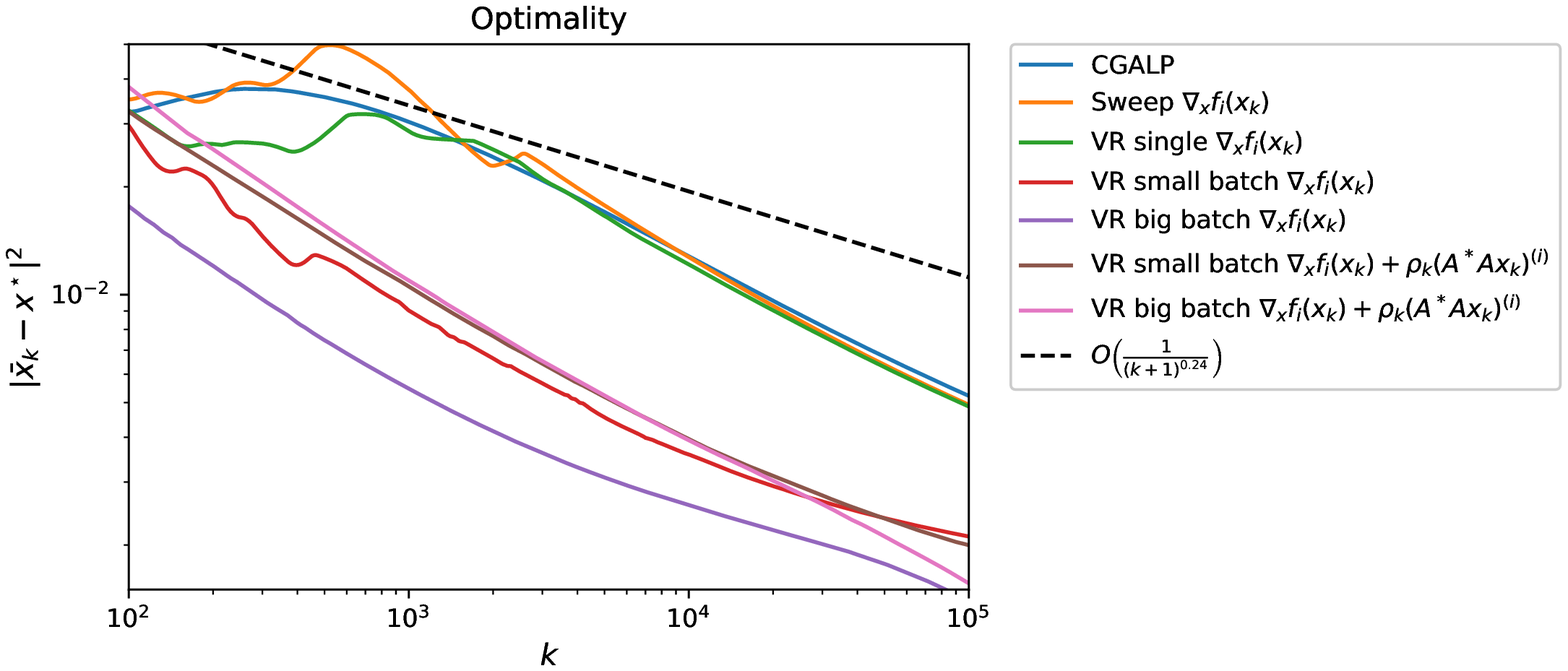}\\
\includegraphics[width=\linewidth]{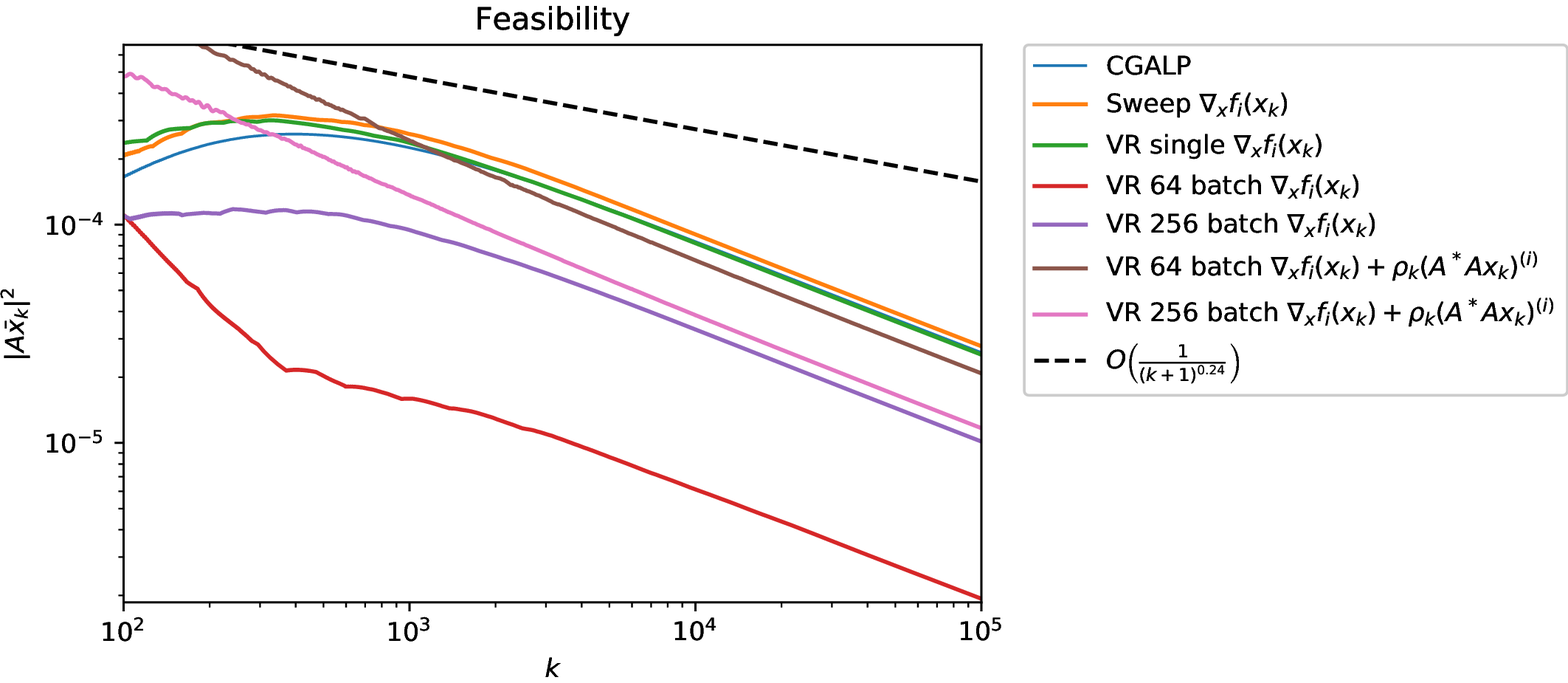}
\caption{Ergodic convergence profiles for \icgalp applied to the projection problem \eqref{numprob} with $n=1024$. The step size is\fekn $\gamma_k = \para{k+1}^{-\para{1-\frac{1}{4} + 0.01}}$ and the weight for variance reduction is\fekn $\nu_k=\gamma_k^{2/3}$.}
\label{fig:bigfig}
\end{figure}

\begin{figure}[hbp]
\centering
\includegraphics[width=\linewidth]{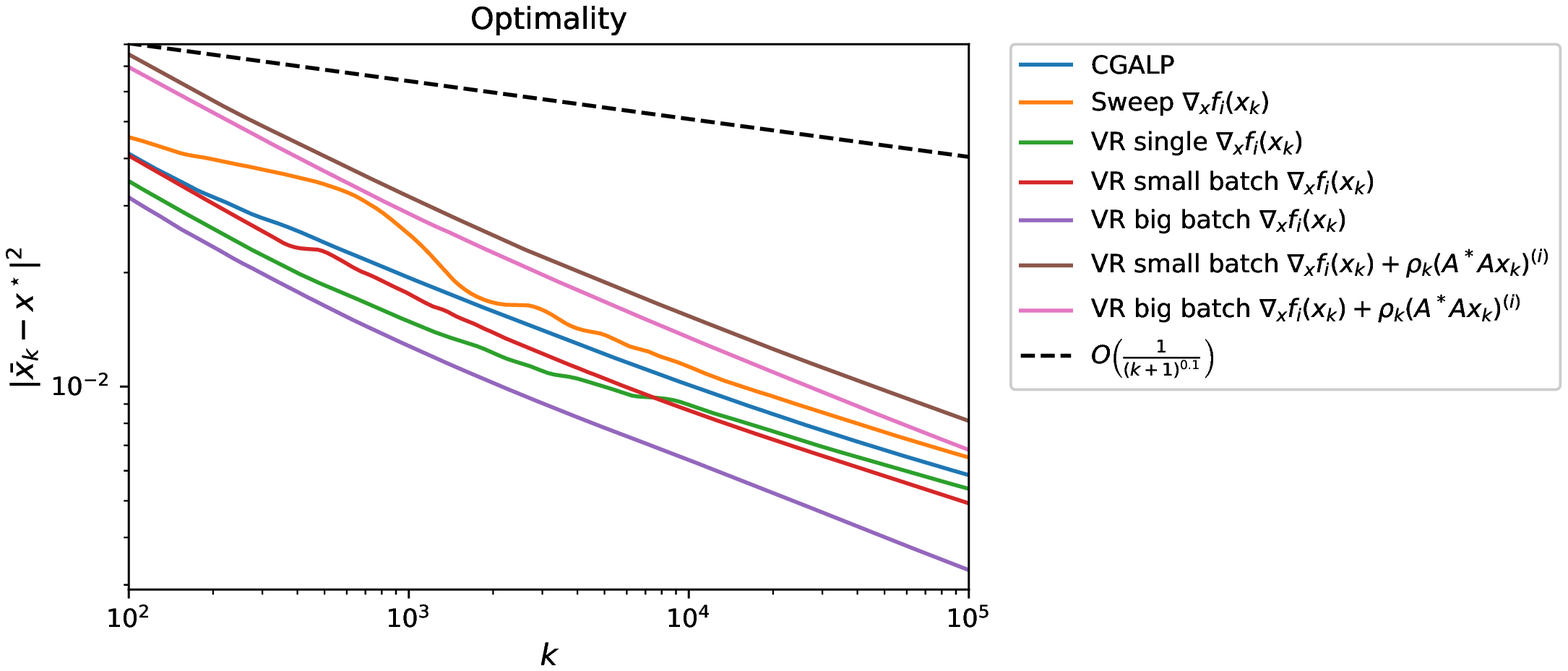}\\
\includegraphics[width=\linewidth]{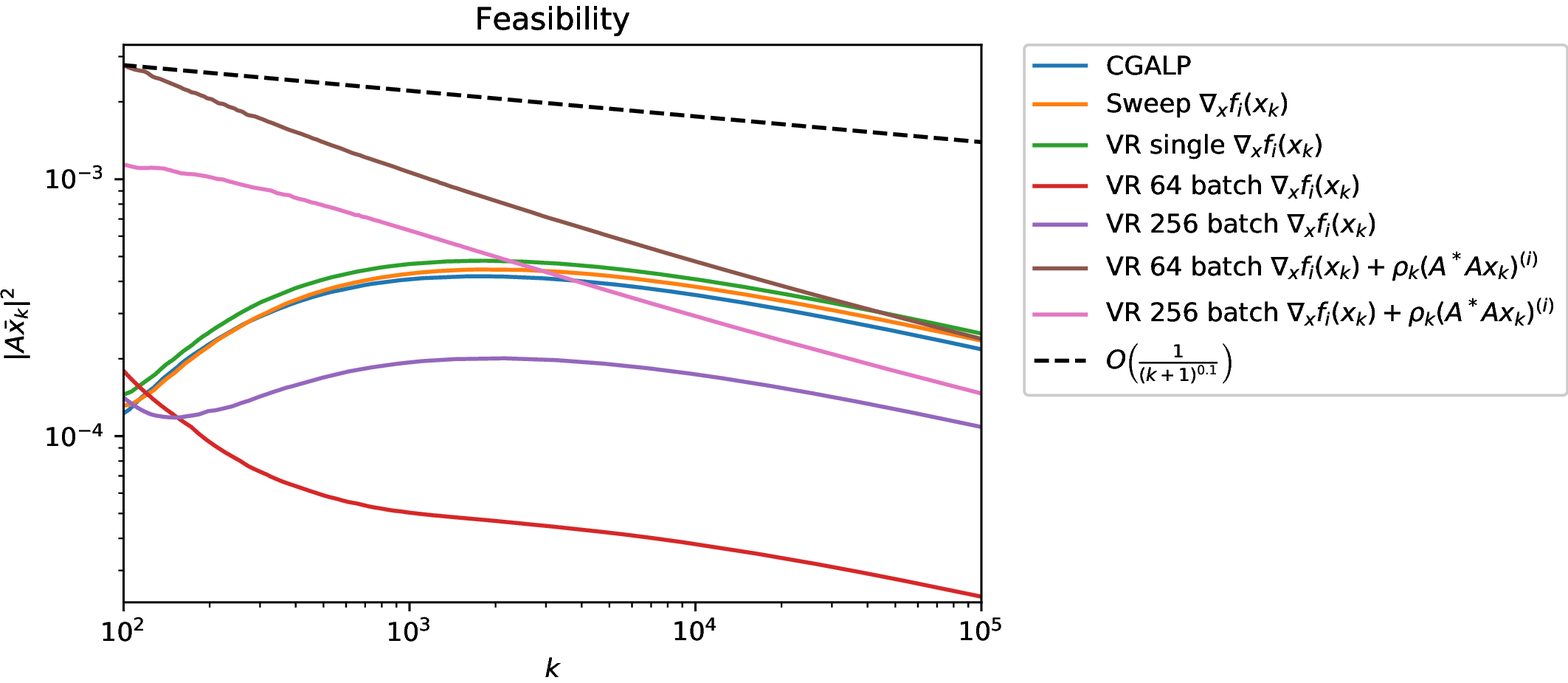}
\caption{Ergodic convergence profiles for \icgalp applied to the projection problem \eqref{numprob} with $n=1024$. The step size is\fekn$\gamma_k = \para{k+1}^{-\para{1-\frac{1}{4}+0.15}}$ and the weight for variance reduction is\fekn $\nu_k=\gamma_k^{2/3}$.}
\label{fig:smallfig}
\end{figure}

Since the problem \eqref{numprob} is strongly convex, we show $\norm{\bar{x}_k-x^\star}{}^2$ in addition to the feasibility gap, $\norm{A\bar{x}_k}{}^2$ where $\bar{x}_k$ is the ergodic variable\fekn
\newq{
\avx_k \eqdef \sum_{i=0}^{k} \gamma_ix_{i+1}/\Gamma_k.
}
We initialize $y\in\R^n$ and $A\in\R^{2\times n}$ randomly. To find the solution $x^\star$ to high precision, we use generalized forward-backward before running the experiments. As a baseline, we run CGALP, the exact counterpart to \icgalp, and display the results. We run the sweeping method on $\nabla f\para{x_k}$ for two different step size choices, displayed in Figures \ref{fig:bigfig} and \ref{fig:smallfig}. For the variance reduction, we examine both the case where $\nabla L\para{x_k,\eta_k}$ is sampled and the case including the gradient of the quadratic term is sampled (see \remref{Esamp}), for two different step size and weight choices as well as different batch sizes ($1, 64,$ or $256$), displayed in Figures \ref{fig:bigfig} and \ref{fig:smallfig}.

\ifdefined\COMPLETE
\else
\documentclass[12pt]{article}
\input{package_header} 
\begin{document}
\fi
\section{Conclusion}\label{sec:conc}
We introduced an inexact extension of the CGALP algorithm, given in \cite{silveti}, which allows for either stochastic or deterministic errors in the computation of several important quantities. The main benefit of this extension will be in the high-dimensional setting, where computing the terms $\nabla f$, $\prox_{\beta g}$, or the linear minimization oracle can be impractical. Several different methods were considered which demonstrated how the gradient $\nabla f$ could be computed in such a way that the summability conditions of \icgalp would be satisfied. The main drawbacks of using the inexact variant of the algorithm emerge from the restrictions on the parameters one is free to choose. Indeed, here the choices of step sizes are more strict than in the CGALP setting. However, the predicted convergence rates for both the optimality and feasibility maintain the same dependence on parameters as was observed for CGALP in an almost sure sense.
\ifdefined\COMPLETE
\else
\end{document}
\fi
\ifdefined\COMPLETE
\else
\documentclass[12pt]{article}
\input{package_header} 
\begin{document}
\fi


\section{Appendix}
\begin{lemma}\label{lem:gamest}
Consider a positive sequence $\seq{u_k}$ which satisfies\fekn
\nnewq{\label{eq:uass}
u_{k+1}\leq \para{1-c\gamma_k^s}u_k + d\gamma_k^t,
}
for some real numbers $s$ and $t$ satisfying $0<s<\min\brac{1,t}$. If, in addition, the sequence $\seq{\gamma_k}$ satisfies\fekn
\nnewq{\label{eq:gamass}
\frac{\gamma_k}{\gamma_{k+1}}\leq 1+o\para{\gamma_k^s},
}
then, for $k$ sufficiently large, it holds,
\newq{
u_k \leq \frac{d}{c} \gamma_k^{t-s}+o\para{\gamma_k^{t-s}}
}
\end{lemma}
\begin{proof}
For each $k\in\N$, we denote $\nu_k\eqdef \gamma_k^{s-t}u_k - \frac{d}{c}$ such that $u_k = \gamma_k^{t-s}\para{\nu_k+\frac{d}{c}}$. Then, by \eqref{eq:uass},
\newq{
\nu_{k+1} = \gamma_{k+1}^{s-t}u_{k+1}-\frac{d}{c} \leq \gamma_{k+1}^{s-t}\para{\para{1-c\gamma_k^s}u_{k} + d\gamma_k^t}-\frac{d}{c} = \gamma_k^{s-t}\para{\frac{\gamma_k}{\gamma_{k+1}}}^{t-s}\para{\para{1-c\gamma_k^s}u_k + d\gamma_k^t}-\frac{d}{c}.
}
By \eqref{eq:gamass}, we then have\fekn
\newq{
\nu_{k+1}\leq \gamma_k^{s-t}\para{1+o\para{\gamma_{k}^s}}^{t-s}\para{\para{1-c\gamma_k^s}u_k+d\gamma_k^t}-\frac{d}{c}.
}
Substituting for $u_k$ using the definition of $\nu_k$ we find\fekn
\newq{
\nu_{k+1}\leq \gamma_k^{s-t}\para{1+o\para{\gamma_k^s}}^{t-s}\para{\para{1-c\gamma_k^s}\para{\nu_k+\frac{d}{c}}\gamma_k^{t-s}+d\gamma_k^t}-\frac{d}{c}.
}
Now, we take a Taylor expansion for the term $\para{1+o\para{\gamma_k}^s}^{t-s}\approx \para{1+o\para{\gamma_k^s}}$ to get, for $k$ sufficiently large,
\newq{
\nu_{k+1}\leq \gamma_k^{s-t}\para{1+o\para{\gamma_k^s}}\para{\para{1-c\gamma_k^s}\para{\nu_k+\frac{d}{c}}\gamma_k^{t-s} + d\gamma_k^t} - \frac{d}{c}.
}
We distribute the $\gamma_k^{s-t}$ and then expand parentheses,
\newq{
\nu_{k+1}&\leq \para{1+o\para{\gamma_k^s}}\para{\para{1-c\gamma_k^s}\para{\nu_k+\frac{d}{c}} + d\gamma_k^s}-\frac{d}{c}\\
&=\para{1-c\gamma_k^s}\nu_k + \para{1-c\gamma_k^s}\frac{d}{c} + d\gamma_k^s + o\para{\gamma_k^s}\para{\para{1-c\gamma_k^s}\para{\nu_k+\frac{d}{c}} + d\gamma_k^s} - \frac{d}{c}\\
&=\para{1-c\gamma_k^s}\nu_k + \para{1-c\gamma_k^s}\frac{d}{c} + d\gamma_k^s + o\para{\gamma_k^s}\para{1-c\gamma_k^s}\nu_k+ o\para{\gamma_k^s}\para{1-c\gamma_k^s}\frac{d}{c} + o\para{\gamma_k^s}d\gamma_k^s - \frac{d}{c} \\
&=\para{1-c\gamma_k^s+ o\para{\gamma_k^s}}\nu_k +o\para{\gamma_k^s}.
}
Fix $0<\tilde c<c$. Then, by definition of $o\para{\gamma_k^s}$, $\exists k_0\in\N$ such that, $\forall k>k_0$, $o\para{\gamma_k^s}\leq (c-\tilde{c})\gamma_k^s$. Then,
\newq{
\para{1-c\gamma_k^s+ o\para{\gamma_k^s}}\nu_k \leq \para{1-\tilde c \gamma_k^s}\nu_k.
}
From this we conclude, by \cite[Ch.2, Lemma~3]{polyak1987introduction}, that $\limsup\limits_k \nu_k \leq 0$.
Thus, by definition of $\nu_{k}$,
\newq{
u_{k+1}\leq \frac{d}{c}\gamma_k^{t-s}+o\para{\gamma_k^{t-s}}.
}
\end{proof}

\ifdefined\COMPLETE
\else
\end{document}
\fi

\section*{Acknowledgements}
ASF was supported by the ERC Consolidated grant NORIA. JF was partly supported by Institut Universitaire de France. CM was supported by Project MONOMADS funded by Conseil R\'egional de Normandie. ASF would like to thank Jingwei Liang for useful discussions had during his visit to Cambridge University.

\appendix
\small
\begin{small}
\bibliographystyle{plain}
\bibliography{cgalp_inexact}
\end{small}
\end{document}